\documentclass[12pt]{amsart}

\usepackage[margin=1.0in]{geometry}
\usepackage[english]{babel}
\usepackage[utf8]{inputenc}
\usepackage{fancyhdr}
\usepackage{natbib}
\usepackage{amsmath}

\RequirePackage{amsthm,amsmath,amsfonts,amssymb}
\usepackage{graphicx}
\usepackage{hyperref}
\usepackage{graphicx}
\usepackage{placeins}
\usepackage{caption}
\usepackage{subcaption}
\usepackage{booktabs}
\numberwithin{equation}{section}

\usepackage{setspace}

\theoremstyle{plain}
\theoremstyle{remark} 
\numberwithin{equation}{section} 

\title{Isotonic Regression Estimators For Simultaneous Estimation  of Order Restricted Location/Scale Parameters of a Bivariate Distribution: A Unified Study}

\author{Naresh Garg  and Neeraj Misra \\ {\footnotesize Department of Mathematics and Statistics\\Indian Institute of Technology Kanpur \\Kanpur-208016, Uttar Pradesh, India}}

\makeatletter
\def\@seccntformat#1{%
	\protect\textup{\protect\@secnumfont
		\ifnum\pdfstrcmp{subsection}{#1}=0 \bfseries\fi
		\csname the#1\endcsname
		\protect\@secnumpunct
	}%
}  
\makeatother

\begin{document}
	\maketitle
	\section*{\textbf{Abstract}}
	
	The problem of simultaneous estimation of location/scale parameters $\theta_1$ and $\theta_2$ of a general bivariate location/scale model, when the ordering between the parameters is known apriori (say, $\theta_1\leq \theta_2$), has been considered. We consider isotonic regression estimators based on the best location/scale equivariant estimators (BLEEs/BSEEs) of $\theta_1$ and $\theta_2$ with general weight functions. Let $\mathcal{D}$ denote the corresponding class of isotonic regression estimators of $(\theta_1,\theta_2)$. Under the sum of the weighted squared error loss function, we characterize admissible estimators within the class $\mathcal{D}$, and identify estimators that dominate the BLEE/BSEE of ($\theta_1$,$\theta_2$). Our study unifies several studies reported in the literature for specific probability distributions having independent marginals. We also report a generalized version of the Katz (1963) result on the inadmissibility of certain estimators under a loss function that is weighted sum of general loss functions for component problems. A simulation study is also carried out to validate the findings of the paper.
	\\~\\ \textbf{Keywords:} Admissible estimators; BLEE; BSEE; Inadmissible estimators; Isotonic regression; Location/Scale model; Mixed estimators.
	
	

	\section{\textbf{Introduction}}

	Let $\underline{X}=(X_1,X_2)$ be a random vector having a joint probability density function (p.d.f.) 
	$f_{\underline{\theta}}(x_1,x_2),\;(x_1,x_2)\in \Re^2, $
	where $\underline{\theta}=(\theta_1,\theta_2)\in \Theta=\Re^2/\Re_{++}^2$ is the vector of unknown location/scale parameters, $\Re_{++}=(0,\infty)$ and $\Re^2$ denotes the two-dimensional Euclidean space. In many real life situations, the ordering between the parameters $\theta_1$ and $\theta_2$ may be known apriori (say, $\theta_1\leq\theta_2$) and it may be of interest to simultaneously estimate $\theta_1$ and $\theta_2$. For example, in an engine efficiency measurement experiment where estimating the average efficiency of an internal combustion engine (I.C. engine) and an external combustion engine (E.C. engine) is of interest, it can be assumed that the average efficiency of an I.C. engine is higher than the average efficiency of an E.C. engine. For more real-life examples, one can see Barlow et al. (\citeyear{MR0326887}) and Robertson et al. (\citeyear{MR961262}).\vspace{3mm}

	Let $\Theta_0=\{\underline{\theta}\in\Theta: \theta_1\leq\theta_2\} $ be the restricted parameter space. There is an extensive literature on estimation of ($\theta_1$,$\theta_2$) (simultaneously as well as componentwise) when it is known apriori that $\underline{\theta}\in \Theta_0$. Most of the early studies were focussed around finding restricted maximum likelihood estimators (under restriction that $\underline{\theta}\in \Theta_0$) and studying their properties (see, for example, van Eeden (\citeyear{MR0083859}, \citeyear{MR0090197}, \citeyear{MR0090196}, \citeyear{MR0102874}) and Brunk (\citeyear{MR73894})). Subsequently, a vast literature appeared on decision theoretic estimation of order restricted estimators with special focus on risk properties of isotonic regression and maximum likelihood estimators (see  Katz (\citeyear{MR150874}), Cohen and Sackrowitz (\citeyear{MR270483}), Brewster-Zidek (\citeyear{MR381098}), Lee (\citeyear{MR615447}), Kumar and Sharma (\citeyear{MR981031}, \citeyear{MR1058934},\citeyear{MR1165706}), Kelly (\citeyear{MR994278}), Kushary and Cohen (\citeyear{MR1029476}), Kaur and Singh (\citeyear{MR1128873}), Gupta and Singh (\citeyear{gupta1992}), Vijayasree and Singh (\citeyear{MR1220404}), Hwang and Peddada (\citeyear{MR1272076}), Kubokawa and Saleh (\citeyear{MR1370413}), Vijayasree et al. (\citeyear{MR1345425}), Misra and Dhariyal (\citeyear{MR1326266}), Garren (\citeyear{MR1802627}), Misra et al. (\citeyear{MR1904424}, \citeyear{MR2205815}), Peddada et al. (\citeyear{MR2202656}), Chang and Shinozaki (\citeyear{chang2015}) and Patra and Kumar (\citeyear{patra2017estimating})). For bivariate models, isotonic regression estimators are also popularly called mixed estimators (see, for example, Katz (\citeyear{MR150874}), Kumar and Sharma (\citeyear{MR981031}), and Vijayasree and Singh (\citeyear{vijayasree1991}, \citeyear{MR1220404})). Most of the above mentioned studies are centred around specific distributions with independent marginals and specific loss functions, barring a few general studies (see, for example, Kelly (\citeyear{MR994278}), Hwang and Peddada (\citeyear{MR1272076}), and Kubokawa and Saleh (\citeyear{MR1370413})). For a detailed account of contributions in this area of research one may refer to the research monograph by van Eeden (\citeyear{MR2265239}).
	\vspace*{3mm}

	A natural goal in these problems is to use the prior information $\underline{\theta}\in \Theta_0$ to find estimators that improve upon the natural estimators for the unrestricted case (e.g., the best location/scale estimators; unrestricted maximum likelihood estimators). This aspect of the problem has been taken up by many researchers in the part for specific probability models and/or specific loss functions (see, for example, Kumar and Sharma (\citeyear{MR1058934}), Kushary and Cohen (\citeyear{MR1029476}), and Vijayasree and Singh (\citeyear{vijayasree1991}, \citeyear{MR1220404})).\vspace*{2mm}

	Katz (\citeyear{MR150874}) introduced mixed estimators for simultaneously estimating order restricted parameters of two distributions. Under a general class of reasonable loss functions (which includes squared error loss), he showed that the natural estimators for the unrestricted case, that are not ordered in the same way as parameters, are inadmissible and he proposed dominating estimators. He further investigated a class of mixed estimators that retain the ordering of parameters and obtained admissible and minimax estimators under normality of the two distributions. Kumar and Sharma (\citeyear{MR981031}) dealt with simultaneous estimation of ordered means of two normal distributions having known variances and characterized admissible estimators among mixed estimators based on the BLEEs of $\theta_1$ and $\theta_2$. Jin and Pal (\citeyear{MR1144879}) considered mixed estimators for simultaneous estimation of order restricted location parameters of two independent exponential distributions and derived estimators that are admissible within the class of mixed estimators. For simultaneous estimation of order restricted scale parameters of two independent exponential distributions, Vijayasree and Singh (\citeyear{vijayasree1991}) considered the class of isotonic regression estimators based on unrestricted maximum likelihood estimators and obtained estimators that are admissible within this class of estimators. Patra and Kumar (\citeyear{patra2017estimating}) considered simultaneous estimation of ordered means of a bivariate normal distribution under the sum of squared error loss functions and discussed admissibility of isotonic regression estimators based on BLEEs of $\theta_1$ and $\theta_2$. \vspace*{2mm}
	
	\vspace*{3mm}

	In this paper we aim to unify some of the above studies for specific probability models by considering a general bivariate location/scale model. For simultaneous estimation of $\theta_1$ and $\theta_2$, under the restricted parameter space and the sum of weighted squared error loss functions, we consider the class of isotonic regression estimators based of BLEEs/BSEEs. We will characterize admissible estimators within this class and obtain subclass of estimators that dominate the BLEEs/BSEEs. The rest of the paper is organized as follows. In Section 2, we introduce some useful definitions and results that are used later in the paper. In Section 3, we report a general result on the inadmissibility of certain estimators under a loss function that is weighted sum of general loss functions. Section 4 (Section 5) deals with isotonic regression estimators of order restricted location (scale) parameters.

	\section{\textbf{Some Useful Definitions and Results}}

	We will discuss some properties of log-concave and log-convex functions (see Pecaric et. al. (\citeyear{MR1162312})), as they are relevant to our study. We begin with the definitions of log-concave and log-convex functions.
	\\~\\ \textbf{Definition 2.2} Let $D$ be a convex subset of $\Re^m$, the $m$-dimensional Euclidean space. A function $h:D\rightarrow [0,\infty)$ is said to be log-concave (log-convex) if, for all $\underline{x},\,\underline{y}\in D$ and all $\alpha\in[0,1]$,
	$$h(\alpha \underline{x}+(1-\alpha) \underline{y})\geq (\leq) (h(\underline{x}))^{\alpha} (h(\underline{y}))^{1-\alpha}.$$
	The following property is well known.
	\\~\\ \textbf{P1.} Let $g:\Re\rightarrow [0,\infty)$ is a p.d.f. with the interior of support as $(a,b)$ $(-\infty\leq a< b\leq \infty)$. Then $g$ is log-concave (log-convex) on $(a,b)$ if, and only if,
	\begin{equation}\label{eq:2.1}
		g(x_1)\,g(x_2-\delta)\geq (\leq) g(x_1-\delta) \,g(x_2),
	\end{equation}
	for all $0<\delta <b-a$ and for all $a+\delta<x_1<x_2<b$ (see Pecaric et. al. (\citeyear{MR1162312})).
	\vspace*{2mm}

	Throughout the paper, whenever, we say that a function $ k: A \rightarrow \Re $, where $ A \subseteq \Re $, is increasing (decreasing) it means that it is non-decreasing (non-increasing). Moreover, $\Re$ will denote the real line, $\Re_{++}^2=(0,\infty)\times (0,\infty)$ and, for any integer $m\;(\geq \,2)$, $\Re^m$ will denote the $m$-dimensional Euclidean space.
	\vspace*{2mm}
	
	The following result, taken from Misra and van der Meulen (\citeyear{MR1985890}), will be used in proving the main results of the paper.
	\\~\\ \textbf{Proposition 2.1} Let $T_1$ and $T_2$ be random variables having distributional supports $B_1$ and $B_2$, respectively. Let $k_i:B_1\cup B_2\rightarrow \Re,\;i=1,2,$ be given functions and let $\beta_i=E[k_i(T_i)],\;i=1,2$. Suppose that $T_1\leq_{st} T_2$ ($T_2\leq_{st} T_1$). Then
	\\ \textbf{(i)} $\beta_1\leq \beta_2$, provided $k_1(x)\leq k_2(x)$, $\forall\; x\in B_1\cup B_2$, and $k_1(x)$ or $k_2(x)$ is an increasing (decreasing) function of $x\in B_1\cup B_2$;
	\\ \textbf{(ii)} $\beta_2\leq \beta_1$, provided $k_2(x)\leq k_1(x)$, $\forall\; x\in B_1\cup B_2$, and $k_1(x)$ or $k_2(x)$ is an decreasing (increasing) function of $x\in B_1\cup B_2$.
	\vspace*{2mm}

	In the following section we consider isotonic regression estimators of order restricted location parameters $\theta_1$ and $\theta_2$ ($\theta_1\leq \theta_2$) under a quite general loss function.
	\vspace*{2mm}

	\section{\textbf{A General Result for Inadmissibility of an Unrestricted Estimator}}
	
	Let $\underline{X}=(X_1,X_2)$ be a random vector having the Lebesgue probability density function (p.d.f.) $f_{\underline{\theta}}(x_1,x_2),\; (x_1,x_2)\in \Re^2, \;\underline{\theta}=(\theta_1,\theta_2)\in \Theta_0,$
	where $\underline{\theta}=(\theta_1,\theta_2)\in\Theta_0=\{(x_1,x_2):-\infty<x_1\leq x_2<\infty\}$ is an unknown parameter. Consider simultaneous estimation of $\theta_1$ and $\theta_2$ under the loss function
	\begin{equation}\label{eq:3.1}
		L(\underline{\theta},\underline{a})=p_1 W(a_1-\theta_1) +p_2 W(a_2-\theta_2),\;\underline{\theta}=(\theta_1,\theta_2)\in\Theta_0,\;\underline{a}=(a_1,a_2)\in\mathcal{A}=\Re^2,
	\end{equation}
	where $p_1>0$ and $p_2>0$ are pre-specified constants and $W:\Re\rightarrow [0,\infty)$ satisfies the following assumption:
	\\~\\\textbf{Assumption D1.}  $W(\cdot)$ is a non-negative and strictly convex function, defined on $\Re$, such that $W(0)=0$, $W(t)$ is strictly decreasing in $t\in (-\infty,0)$ and strictly increasing in $t\in (0,\infty)$.
	\vspace*{2mm}
	
	Under the assumption $D1$, $W(\cdot)$ is continuous everywhere and it is also differentiable everywhere, except possibly at countable number of points. 
	\vspace*{2mm}
	
	\noindent For an estimator $\underline{\delta}=(\delta_1,\delta_2)$, its risk function is defined by 
	$$R(\underline{\theta},\underline{\delta})=E_{\underline{\theta}}[L(\underline{\theta},\underline{\delta})]=p_1 E_{\underline{\theta}}[W_1(\delta_1(\underline{X})-\theta_1)]+ p_2 E_{\underline{\theta}}[W_2(\delta_2(\underline{X})-\theta_2)],\; \; \underline{\theta}\in \Theta_0.$$
	\vspace*{1.5mm}
	
	\noindent The following lemma generalizes the result stated in Theorem 1 of Katz (1963).
	
	\vspace*{3mm}
	
	\noindent \textbf{Lemma 3.1.} Suppose that the assumption (D1) holds. Let $\underline{\delta}=(\delta_1,\delta_2)\in \Re^2$ be such that $\delta_1>\delta_2$. For $\alpha\in \Re$, define $\delta_{1,\alpha}=\alpha \delta_1 +(1-\alpha)\delta_2$, $\delta_{2,\alpha}=\frac{p_1}{p_2}(1-\alpha) \delta_{1}+(1-\frac{p_1}{p_2}(1-\alpha))\delta_{2}$ and $\underline{\delta}_{\alpha}=(\delta_{1,\alpha},\delta_{2,\alpha})$. Then,
	\\~\\ \textbf{(a)} for $\max\{0,\frac{p_1-p_2}{p_1}\}< \alpha<1$, $$L(\underline{\theta},\underline{\delta}_{\alpha})< L(\underline{\theta},\underline{\delta}),\; \forall \; \underline{\theta}\in \Theta_0;$$
	\textbf{(b)} for $\max\{0,\frac{p_1-p_2}{p_1}\}\leq \alpha<1$, $$L(\underline{\theta},\underline{\delta}_{\alpha})< L(\underline{\theta},\underline{\delta}),\; \forall \; \underline{\theta}\in \Theta_0, \text{ such that } \theta_1<\theta_2;$$
	
	\noindent \textbf{(c)} for $\frac{p_1}{p_1+p_2} \leq \alpha_1 <\alpha_2<\infty$,
	$$L(\underline{\theta},\underline{\delta}_{\alpha_1})< L(\underline{\theta},\underline{\delta}_{\alpha_2}),\; \forall \; \underline{\theta}\in \Theta_0.$$

	\noindent
	\textbf{Proof.} \textbf{(a) and (b)}  Since $W(\cdot)$ is a strictly convex function $\Re$, we have
	\begin{equation}\label{eq:3.2}
		\frac{W(x_2)-W(x_1)}{x_2-x_1}<    \frac{W(y_2)-W(y_1)}{y_2-y_1}	,
	\end{equation}
	provided $x_1\leq y_1,\; x_2\leq y_2$, with at least one of the inequalities being strict, $x_1\neq x_2$ and $y_1\neq y_2$. Fix $\alpha\in [\max\{0,\frac{p_1-p_2}{p_1}\},1)$, and $\underline{\theta}\in \Theta_0$. Using (3.2), with $x_1=\delta_2-\theta_2$, $x_2=\frac{p_1}{p_2}(1-\alpha) \delta_1+(1-\frac{p_1}{p_2}(1-\alpha))\delta_2-\theta_2$, $y_1=\alpha \delta_1+(1-\alpha)\delta_2-\theta_1$ and $y_2=\delta_1-\theta_1$ (so that $x_1<y_1$, $x_2<y_2$, $x_1\neq x_2$ and $y_1\neq y_2$), we get
	
	\begin{align*}
		L(\underline{\theta},\underline{\delta})-L(\underline{\theta},\underline{\delta}_{\alpha})	=& p_1 [W(\delta_1-\theta_1)-W(\alpha \delta_1+(1-\alpha)\delta_2-\theta_1)]\\
		& +p_2\left[W(\delta_2-\theta_2)-W\left(\frac{p_1}{p_2}(1-\alpha) \delta_1+(1-\frac{p_1}{p_2}(1-\alpha))\delta_2-\theta_2\right)\right]\\
		=&
		p_1(1-\alpha)(\delta_1-\delta_2) \bigg\{	\frac{W(\delta_1-\theta_1)-W(\alpha \delta_1+(1-\alpha)\delta_2-\theta_1)}{(1-\alpha)(\delta_1-\delta_2)}\\
		&\; - \frac{W(\frac{p_1}{p_2}(1-\alpha) \delta_1+(1-\frac{p_1}{p_2}(1-\alpha))\delta_2-\theta_2)-W(\delta_2-\theta_2)}{(1-\alpha)\frac{p_1}{p_2}(\delta_1-\delta_2)} \bigg\}\\
		>&\; 0,
	\end{align*}
	under the hypotheses of assertions (a) and (b).\vspace*{2mm}
	
	\noindent \textbf{(c)}  Fix $\underline{\theta}\in \Theta_0$ and consider
	\begin{align*}
		Q(\alpha)= L(\underline{\theta},\underline{\delta}_{\alpha})= &p_1 [W(\alpha \delta_1+(1-\alpha)\delta_2-\theta_1)]\\&+ p_2 \left[W\left(\frac{p_1}{p_2}(1-\alpha) \delta_1+(1-\frac{p_1}{p_2}(1-\alpha))\delta_2-\theta_2\right)\right].
	\end{align*}
	For simplicity, we assume that $W(\cdot)$ is differentiable everywhere, so that $W^{'}(\cdot)$ is strictly increasing (as $W(\cdot)$ is convex), Then
	\begin{align*}
		Q^{'}(\alpha)=p_1(\delta_1-\delta_2) &\bigg\{W^{'}(\alpha \delta_1+(1-\alpha)\delta_2-\theta_1)\\&\qquad\qquad - W^{'}\left(\frac{p_1}{p_2}(1-\alpha) \delta_1+(1-\frac{p_1}{p_2}(1-\alpha))\delta_2-\theta_2\right)\bigg\}>0,
	\end{align*}
	provided $\alpha\geq \frac{p_1}{p_1+p_2}$ and $\theta_1\leq \theta_2$, with at least one of the inequalities being strict; here $Q^{'}(\cdot)$ denotes the derivative of $Q(\cdot)$. \vspace*{2mm}

	\noindent The following theorem is an immediate consequence of the above lemma.
	\vspace*{4mm}

	\noindent \textbf{Theorem 3.1.} Suppose that the assumption (D1) holds. Let $\underline{\delta}(\underline{X})=(\delta_1(\underline{X}),\delta_2(\underline{X}))$ be an estimator of $\underline{\theta}=(\theta_1,\theta_2)$ such that $P_{\underline{\theta}}[\delta_1(\underline{X})>\delta_2(\underline{X})]>0$, for some $\underline{\theta}\in \Theta_0$, such that $\theta_1<\theta_2$. For $\alpha\in \Re$, define $\underline{\delta}_{\alpha}^{M}=(\delta_{1,\alpha}^{M},\delta_{2,\alpha}^{M}),$ where 
	
	\begin{align*} \delta_{1,\alpha}^{M}(\underline{X})=\begin{cases}
			\delta_{1}(\underline{X}), &\text{ if}\quad \delta_{1}(\underline{X})\leq \delta_{2}(\underline{X}) \\
			\alpha \delta_{1}(\underline{X})+(1-\alpha)\delta_{2}(\underline{X}),&\text{ if}\quad \delta_{1}(\underline{X})> \delta_{2}(\underline{X})
		\end{cases} 
	\end{align*}
	and                                                                                    \begin{align*}	\delta_{2,\alpha}^{M}(\underline{X})=\begin{cases}
			\delta_{2}(\underline{X}), &\text{ if}\quad \delta_{1}(\underline{X})\leq \delta_{2}(\underline{X}) \\
			\frac{p_1}{p_2}(1-\alpha) \delta_{1}(\underline{X})+(1-\frac{p_1}{p_2}(1-\alpha))\delta_{2}(\underline{X}),&\text{ if}\quad \delta_{1}(\underline{X})> \delta_{2}(\underline{X})
		\end{cases} .
	\end{align*}     
	Then,
	\\~\\ \textbf{(a)} for $\max\{0,\frac{p_1-p_2}{p_1}\}\leq \alpha< 1$, $$R(\underline{\theta},\underline{\delta}_{\alpha}^M)\leq R(\underline{\theta},\underline{\delta}),\; \forall \; \underline{\theta}\in \Theta_0, $$
	with strict inequality for some $\underline{\theta}\in \Theta_0$;
	\vspace*{2mm}
	
	\noindent \textbf{(b)} for $\frac{p_1}{p_1+p_2} \leq \alpha_1 <\alpha_2<\infty$,
	$$R(\underline{\theta},\underline{\delta}_{\alpha_1}^M)\leq R(\underline{\theta},\underline{\delta}_{\alpha_2}^M),\; \forall \; \underline{\theta}\in \Theta_0,$$
	with strict inequality for some $\underline{\theta}\in \Theta_0$.
	\vspace*{2mm}
	
	\noindent
	\textbf{Proof.}  Note that, for any $\alpha\in \Re$, $L(\underline{\theta},\underline{\delta}_{\alpha}^M(\underline{X}))= L(\underline{\theta},\underline{\delta}(\underline{X})),\; \forall \; \underline{\theta}\in \Theta_0,$ if $\delta_1(\underline{X})\leq \delta_2(\underline{X})$. Now the result follows on using Lemma 3.1.
	\vspace*{4mm}

	Let $\underline{\delta}_{\alpha}^{M}=(\delta_{1,\alpha}^{M},\delta_{2,\alpha}^{M}),\; \alpha\in \Re,$ be as defined in Theorem 3.1. We call $\underline{\delta}_{\alpha}^{M},\; \alpha\in \Re,$ a mixed estimator (or isotonic regression estimator) of $\underline{\theta}=(\theta_1,\theta_2)$ based on $\underline{\delta}=(\delta_1,\delta_2)$. The following points are noteworthy.
	\vspace*{2mm}
	
	\noindent \textbf{Remark 3.1 (a)} We have $\underline{\delta}_{1}^{M}=(\delta_{1,1}^{M},\delta_{2,1}^{M})=(\delta_1,\delta_2).$ Also, for $\alpha\leq 1$, $ \delta_{1,\alpha}^{M}(\underline{X})=\min\{
	\delta_{1}(\underline{X}),
	\alpha \delta_{1}(\underline{X})+(1-\alpha)\delta_{2}(\underline{X})\}$
	and                                                                                   	$\delta_{2,\alpha}^{M}(\underline{X})=\max\{
	\delta_{2}(\underline{X}),
	\frac{p_1}{p_2}(1-\alpha) \delta_{1}(\underline{X})+(1-\frac{p_1}{p_2}(1-\alpha))\delta_{2}(\underline{X})\}$.
	\\~\\ \textbf{(b)}  For $\alpha > \frac{p_1}{p_1+p_2}=\alpha_0$, say, $P_{\underline{\theta}}[\delta_{1,\alpha}^{M}(\underline{X})>\delta_{2,\alpha}^{M}(\underline{X})]=P_{\underline{\theta}}[\delta_{1}(\underline{X})>\delta_2(\underline{X})],\; \forall \; \underline{\theta}\in \Theta_0$. Thus any estimator $\underline{\delta}_{\alpha}^{M}=(\delta_{1,\alpha}^{M},\delta_{2,\alpha}^{M})$, with $\alpha>\frac{p_1}{p_1+p_2}$, is inadmissible for estimating $\underline{\theta}\in \Theta_0$, provided $P_{\underline{\theta}}[\delta_{1}(\underline{X})>\delta_2(\underline{X})]> 0,$ for some $\underline{\theta}\in \Theta_0$. In this case the estimator $\underline{\delta}_{\alpha_0}^{M}=(\delta_{1,\alpha_0}^{M},\delta_{2,\alpha_0}^{M})$ improves upon $\underline{\delta}_{\alpha}^{M}=(\delta_{1,\alpha}^{M},\delta_{2,\alpha}^{M})$, for any $\alpha>\frac{p_1}{p_1+p_2}$.
	\\~\\ \textbf{(c)} Suppose that the parameter space is $\Theta_0=\{(\theta_1,\theta_2)\in \Re_{++}^2: \theta_1\leq \theta_2\}$. For estimating $\underline{\theta}=(\theta_1,\theta_2)$, consider the loss function
	\begin{equation*}
		L(\underline{\theta},\underline{a})=p_1 W\left(\frac{a_1}{\theta_1}-1\right) +p_2 W\left(\frac{a_2}{\theta_2}-1\right), \;\underline{\theta}=(\theta_1,\theta_2)\in\Theta_0,\;\underline{a}=(a_1,a_2)\in\mathcal{A}=\Re_{++}^2,
	\end{equation*}
	where $p_1>0$ and $p_2>0$ are pre-specified constants and $W:\Re\rightarrow [0,\infty)$ is a strictly convex function, such that  $W(1)=0$, $W(t)$ is strictly decreasing in $t\in (-\infty,1)$ and strictly increasing in $t\in (1,\infty)$. Then Lemma 3.1 holds, provided $\delta_1>\delta_2>0$. Also, Theorem 3.1 holds, provided $P_{\underline{\theta}}[\delta_i(\underline{X})>0]=1,\;\forall \; \underline{\theta}\in \Theta_0,\;i=1,2.$
	\vspace*{3mm}

	\section{\textbf{Admissibility of Isotonic regression estimators of Location Parameters under weighted sum of squared errors loss}}
	
	Let $\underline{X}=(X_1,X_2)$ be a random vector having the Lebesgue probability density function (p.d.f.) belonging to the location family
	\begin{equation}\label{eq:4.1}
		f_{\underline{\theta}}(x_1,x_2)= 	f(x_1-\theta_1,x_2-\theta_2),\; (x_1,x_2)\in \Re^2, \;\underline{\theta}=(\theta_1,\theta_2)\in \Theta_0,
	\end{equation} 
	where $f(\cdot,\cdot) $ is a specified bivariate Lebesgue p.d.f. on $\Re^2=(-\infty,\infty)\times(-\infty,\infty)$ and $\Theta_0=\{(\theta_1,\theta_2):-\infty<\theta_1\leq \theta_2<\infty\}$. Generally, $(X_1,X_2)$ is a minimal sufficient statistic based on a random sample.  For simultaneous estimation of $\theta_1$ and $\theta_2$ ($\underline{\theta}\in \Theta_0$), consider the weighted sum of squared errors loss function, given by
	\begin{equation}\label{eq:4.2}
		L(\underline{\theta},\underline{a})=p_1 (a_1-\theta_1)^2 +p_2 (a_2-\theta_2)^2,\;\underline{\theta}=(\theta_1,\theta_2)\in\Theta_0,\;\underline{a}=(a_1,a_2)\in\mathcal{A}=\Re^2,
	\end{equation}
	where $p_1>0$ and $p_2>0$ are pre-specified constants.
	\vspace*{2mm}

	Under the (unrestricted) parameter space $\Theta=\Re^2$, the problem of simultaneously estimating $\theta_1$ and $\theta_2$ under the loss function (4.2) is invariant under the additive group of transformations $\mathcal{G}=\{g_{c_1,c_2}:(c_1,c_2)\in\Re^2 \}$, where $g_{c_1,c_2}(x_1,x_2)=(x_1+c_1,x_2+c_2)$, $(x_1,x_2)\in\Re^2,\;  (c_1,c_2)\in\Re^2$. The group of transformations $\mathcal{G}$ induces the group of transformations $\overline{\mathcal{G}}=\{\overline{g}_{c_1,c_2}:(c_1,c_2)\in\Re^2 \}$ and $\tilde{\mathcal{G}}=\{\tilde{g}_{c_1,c_2}:(c_1,c_2)\in\Re^2 \}$ on the parameter space $\Theta$ and the action space $\mathcal{A}$, respectively, where $\overline{g}_{c_1,c_2}(x_1,x_2)=(x_1+c_1,x_2+c_2)$, $\tilde{g}_{c_1,c_2}(a_1,a_2)=(a_1+c_1,a_2+c_2)$, $(x_1,x_2)\in\Theta=\Re^2,\;(a_1,a_2)\in\mathcal{A}=\Re^2,\;  (c_1,c_2)\in\Re^2$. Under the group of transformations $\mathcal{G}$, an estimator $\underline{\delta}=(\delta_1,\delta_2)$ is invariant (location invariant) for estimating $\underline{\theta}=(\theta_1,\theta_2)$ if, and only if,
	\begin{equation}\label{eq:4.3}
		\underline{\delta}(\underline{X})=(\delta_1(\underline{X}),\delta_2(\underline{X}))=(X_1-c_1,X_2-c_2),
	\end{equation}
	for some $\underline{c}=(c_1,c_2)\in \Re^2$. The best location equivariant estimator (BLEE) of $(\theta_1,\theta_2)$ is 
	\begin{equation}\label{eq:4.4}
		\underline{\delta}_{0}(\underline{X})=(\delta_{1,0}(\underline{X}),\delta_{2,0}(\underline{X}))=(X_1-c_{0,1},X_2-c_{0,2}),
	\end{equation}
	where, for $Z_i=X_i-\theta_i$, $c_{0,i}=E[Z_i],\;i=1,2$.

	\vspace*{3mm}

	\noindent Consider isotonic regression estimators $\underline{\delta}_{\alpha}^{M}=(\delta_{1,\alpha}^{M},\delta_{2,\alpha}^{M})\; (\alpha\in \Re),$ based on the BLEE $(X_1-c_{0,1},X_2-c_{0,2})$, where, for $\alpha\in \Re$,
	\begin{align}\label{eq:4.5}
		\delta_{1,\alpha}^{M}(\underline{X})=\begin{cases}
			X_1-c_{0,1}, &\text{ if}\quad X_1-c_{0,1}\leq X_2-c_{0,2} \\
			\alpha (X_1-c_{0,1})+(1-\alpha)(X_2-c_{0,2}),&\text{ if}\quad X_1-c_{0,1}> X_2-c_{0,2}
		\end{cases} 
	\end{align}
	and         
	\begin{align}\label{eq:4.6}
		\delta_{2,\alpha}^{M}(\underline{X})=\begin{cases}
			X_2-c_{0,2}, &\text{ if}\quad X_1-c_{0,1}\leq X_2-c_{0,2} \\
			\frac{p_1}{p_2}(1-\alpha) (X_1-c_{0,1})+(1-\frac{p_1}{p_2}(1-\alpha))(X_2-c_{0,2}),&\text{ if}\quad X_1-c_{0,1}> X_2-c_{0,2}
		\end{cases}.
	\end{align}   
	\vspace*{3mm}
	
	Let $\mathcal{D}=\{\underline{\delta}_{\alpha}^{M}=(\delta_{1,\alpha}^{M}(\underline{X}),\delta_{2,\alpha}^{M}(\underline{X})):\alpha\in\Re\}$ be the class of isotonic regression estimators (IREs), based on BLEEs. Define $\lambda=\theta_2-\theta_1$ ($\lambda\geq 0$), $Z_1=X_1-\theta_1$, $Z_2=X_2-\theta_2$, $Z=Z_2-Z_1$ and $\underline{Z}=(Z_1,Z_2)$. The risk function of the estimator $\underline{\delta}_{\alpha}^{M},\; \alpha \in \Re,$ is given by
	\begin{small}
		\\~\\$R(\underline{\theta},\underline{\delta}_{\alpha}^{M})$
		\begin{align}\label{eq:4.7}
			\nonumber
			&=p_1E_{\underline{\theta}}[(Z_1-c_{0,1})^2 I_{(c_{0,2}-c_{0,1}-\lambda,\infty)}(Z)]\\ \nonumber
			&\quad +p_1 E_{\underline{\theta}}[((Z_2+\lambda-c_{0,2})-\alpha (Z+\lambda+c_{0,1}-c_{0,2}))^2 I_{(-\infty,c_{0,2}-c_{0,1}-\lambda)}(Z)]\\ \nonumber
			&\quad +p_2E_{\underline{\theta}}[(Z_2-c_{0,2})^2 I_{(c_{0,2}-c_{0,1}-\lambda,\infty)}(Z)]\\
			&\quad + p_2E_{\underline{\theta}}\left[\!\left(\frac{p_1}{p_2}\alpha (Z+\lambda+c_{0,1}-c_{0,2})-\frac{p_1}{p_2} (Z+\lambda+c_{0,1}-c_{0,2})+Z_2-c_{0,2} \right)^2\! I_{(-\infty,c_{0,2}-c_{0,1}-\lambda)}(Z)\!\right],\, \underline{\theta}\in \Theta_0,
		\end{align}
	\end{small}
	where, for any set $A$, $I_A(\cdot)$ denotes its indicator function.
	The risk function $R(\underline{\theta},\underline{\delta}_{\alpha}^{M})$ depends on $\underline{\theta}\in\Theta_0$ only through $\lambda=\theta_2-\theta_1\;(\geq 0)$. Let $f_Z(\cdot)$ denote the p.d.f. of $Z$, so that $f_Z(z)=\int_{-\infty}^{\infty} f(s,s+z)\,ds,\; -\infty<z<\infty$. 
	Clearly, for any fixed $\underline{\theta}\in\Theta_0$ (or $\lambda\geq 0)$, $R(\underline{\theta},\underline{\delta}_{\alpha}^{M})$ is minimized at $\alpha= \alpha_1(\lambda)$, where, for $\lambda \geq 0$,
	\begin{align}\label{eq:4.8}
		\alpha_1(\lambda)&=\frac{p_1}{p_1+p_2}+\frac{p_2}{p_1+p_2} \frac{\lambda \,E_{\underline{\theta}}[(Z+\lambda+c_{0,1}-c_{0,2}) I_{(-\infty,c_{0,2}-c_{0,1}-\lambda)}(Z)]}{E_{\underline{\theta}}[(Z+\lambda+c_{0,1}-c_{0,2})^2   \nonumber I_{(-\infty,c_{0,2}-c_{0,1}-\lambda)}(Z)]}\\
		&=\frac{p_1}{p_1+p_2}+\frac{p_2}{p_1+p_2} \alpha_1^{*}(\lambda);
	\end{align}
	\begin{align}\label{eq:4.9}
		\text{here}\;\;\;	\alpha_1^{*}(\lambda)&=\frac{\lambda \,E_{\underline{\theta}}[(Z+\lambda+c_{0,1}-c_{0,2}) I_{(-\infty,c_{0,2}-c_{0,1}-\lambda)}(Z)]}{E_{\underline{\theta}}[(Z+\lambda+c_{0,1}-c_{0,2})^2           \nonumber I_{(-\infty,c_{0,2}-c_{0,1}-\lambda)}(Z)]},\;\;\lambda\geq 0. \qquad \qquad\\ \nonumber
		&= \frac{\lambda \,\int_{-\infty}^{0} z\,f_Z(z+c_{0,2}-c_{0,1}-\lambda) dz}{\int_{-\infty}^{0} z^2\,f_Z(z+c_{0,2}-c_{0,1}-\lambda) dz},\;\; \lambda\geq 0,\\
		&=E_{\lambda}[k(S_{\lambda},\lambda)],\;\; \lambda\geq 0,
	\end{align}
	where $k(z,\lambda)=\frac{\lambda}{z},\; z<0,\; \lambda\geq 0,$ and $S_{\lambda}$ is a r.v. having the p.d.f.
	\begin{equation}\label{eq:4.10}
		h_{\lambda}(z)=\begin{cases}
			\frac{z^2 f_Z(z+c_{0,2}-c_{0,1}-\lambda)}{\int_{-\infty}^{0} s^2 f_Z(s+c_{0,2}-c_{0,1}-\lambda)\,ds},& \text{if}\;\; z<0\\
			0, &\text{otherwise}
		\end{cases},\;\;\lambda\geq 0.
	\end{equation}
	
	Using (2.1), it is easy to check that if $f_Z(\cdot)$ is log-concave on $(-\infty,c_{0,2}-c_{0,1})$, then $S_{\lambda_1}\leq_{lr}S_{\lambda_2}$, and consequently $S_{\lambda_1}\leq_{st}S_{\lambda_2}$, whenever $0\leq \lambda_1<\lambda_2<\infty$. \vspace*{2mm}
	
	\noindent
	The following lemma will be useful in proving the main result of this subsection.
	\\~\\ \textbf{Lemma 4.1.} Suppose that $f_Z(\cdot)$ is log-concave on $(-\infty,c_{0,2}-c_{0,1})$. Then $\alpha_1^{*}(\lambda)$ (and hence $\alpha_1(\lambda)$) is a decreasing function of $\lambda\in[0,\infty)$,
	$\inf_{\lambda\geq 0} \alpha_1(\lambda) = \frac{p_1}{p_1+p_2}+\frac{p_2}{p_1+p_2}\lim_{\lambda \to \infty}\alpha_1^{*}(\lambda)=\alpha_{\infty}, \text{ say, and }
	\sup_{\lambda\geq 0} \alpha_1(\lambda) =\frac{p_1}{p_1+p_2}+\frac{p_2}{p_1+p_2}\alpha_1^{*}(0)=\frac{p_1}{p_1+p_2}.$
	\begin{proof}  Let $0\leq \lambda_1<\lambda_2 <\infty$, $T_i=S_{\lambda_i},\;i=1,2,$ and $m_{i}(z)=k(z,\lambda_i)=\frac{\lambda_i}{z},\;z<0,\;i=1,2,$ where $S_{\lambda}$ is a r.v. having p.d.f. given by \eqref{eq:4.10}. Then $\alpha_1^{*}(\lambda_i)=E_{\lambda_i}[m_{i}(T_i)],\;i=1,2$. The hypothesis that $f_Z(\cdot)$ is log-concave on $(-\infty,c_{0,2}-c_{0,1})$ implies that $T_1\leq_{lr} T_2$ and, consequently, $T_1\leq_{st} T_2$. Also $m_{1}(t) \geq \,m_{2}(t),\;\forall \; t<0$, and $m_{i}(t)$ is a decreasing function of $t\in(-\infty,0)$. Now, using Proposition 2.1 (a), it follows that
		$$\alpha_1^{*}(\lambda_1)=E_{\lambda_1}[m_{1}(T_1)]\geq E_{\lambda_2}[m_{2}(T_2)]=\alpha_1^{*}(\lambda_2).$$
	\end{proof}
	\noindent In the following theorem we use the convention that $[a_{\infty},\frac{p_1}{p_1+p_2}]=(a_{\infty},\frac{p_1}{p_1+p_2}]$, if $a_{\infty}=-\infty$. \vspace*{2mm}

	\noindent  \textbf{Theorem 4.1.}  Suppose that $f_Z(\cdot)$ is log-concave on $(-\infty,c_{0,2}-c_{0,1})$. Then the estimators that are admissible within the class $\mathcal{D}=\{\underline{\delta}_{\alpha}^{M}=(\delta_{1,\alpha}^{M},\delta_{2,\alpha}^{M}): \alpha \in \Re\}$ are $\{\underline{\delta}_{\alpha}^{M}=(\delta_{1,\alpha}^{M},\delta_{2,\alpha}^{M}): \alpha \in [\alpha_{\infty},\frac{p_1}{p_1+p_2}]\}$. Moreover, for $-\infty<\alpha_1<\alpha_2\leq \alpha_{\infty}<\infty$ or $\frac{p_1}{p_1+p_2} \leq \alpha_2 < \alpha_1<\infty$, the estimator $\underline{\delta}_{\alpha_2}^{M}(\underline{X})$ dominates the estimator $\underline{\delta}_{\alpha_1}^{M}(\underline{X})$, for any $\underline{\theta}\in \Theta_0$.
	\begin{proof} Let $\alpha_1(\lambda)$ be as defined by \eqref{eq:4.8}, so that, for any fixed $\underline{\theta}\in\Theta_0$ (or fixed $\lambda\geq 0$), the risk function $R(\underline{\theta},\underline{\delta}_{\alpha}^{M})$, given by \eqref{eq:4.7}, is uniquely minimized at $\alpha=\alpha_1(\lambda)$. Since, for any $\lambda \geq 0$, $\alpha_1(\lambda)$ is a continuous function of $\lambda\in[0,\infty)$, it assumes all values in $\Re$ that are between $\inf_{\lambda\geq 0} \alpha_1(\lambda)=\alpha_{\infty}$ and $\sup_{\lambda\geq 0}\alpha_1(\lambda)=\frac{p_1}{p_1+p_2}$, as $\lambda$ varies on $[0,\infty)$. It follows that each $\alpha\in(\alpha_{\infty},\frac{p_1}{p_1+p_2}]$ uniquely minimizes the risk function $R(\underline{\theta},\underline{\delta}_{\alpha}^{M})$ at some $\underline{\theta}\in\Theta_0$ (or at some $\lambda\geq 0$). This proves that the estimators $\{\underline{\delta}_{\alpha}^{M}:\alpha\in(\alpha_{\infty},\frac{p_1}{p_1+p_2}]\}$ are admissible among the estimators in the class $\mathcal{D}$. When $\alpha_{\infty}>-\infty$, the admissibility of the estimator $\underline{\delta}_{\alpha_{\infty}}^M$, within the class $\mathcal{D}$ of isotonic regression estimators, can be proved by contradiction. Also, note that, for any fixed $\lambda\geq 0$, $R(\underline{\theta},\underline{\delta}_{\alpha}^M)$ is a strictly decreasing function of $\alpha$ on $(-\infty,\alpha_1(\lambda))$ and a strictly increasing function of $\alpha$ on $(\alpha_1(\lambda),\infty)$. Since $\alpha_1(\lambda)\in [\alpha_{\infty},\frac{p_1}{p_1+p_2}],\;\forall\;\lambda\geq 0$, the foregoing discussion implies that, for any $\underline{\theta}\in\Theta_0$ (or for any $\lambda\geq 0$), $R(\underline{\theta},\underline{\delta}_{\alpha}^{M})$ is a decreasing function of $\alpha$ on $(-\infty,\alpha_{\infty}]$ (provided $\alpha_{\infty}>-\infty$) and is an increasing function of $\alpha$ on $[\frac{p_1}{p_1+p_2},\infty)$. This establishes the second assertion.
	\end{proof}

	\vspace*{3mm}
	
	\noindent \textbf{Remark 4.1} Note that $\underline{\delta}_{0}(\underline{X})=(\delta_{1,0}(\underline{X}),\delta_{2,0}(\underline{X}))=(X_1-c_{0,1},X_2-c_{0,2})$ is the BLEE of $\underline{\theta}=(\theta_1,\theta_2)$. If $P_{\underline{\theta}}[\delta_{1,0}(\underline{X})>\delta_{2,0}(\underline{X})]=P_{\underline{\theta}}[X_1-c_{0,1}>X_2-c_{0,2}]>0,$ for some $\underline{\theta}\in \Theta_0$, then the mixed estimator $\underline{\delta}_{\alpha}^{M}=(\delta_{1,\alpha}^{M},\delta_{2,\alpha}^{M})$ improves upon $\underline{\delta}_{0}(\underline{X})$, provided $\max\{0,\frac{p_1-p_2}{p_1}\}\leq \alpha <1$.
	
	\vspace*{3mm}

	\subsection{Applications}
	\label{sec:3.3}
	\noindent
	
	\vspace*{3mm}

	\noindent
	Now we illustrate some applications of results of Sections 3 and 4.
	
	\vspace*{3mm}

	\noindent
	\textbf{Example 4.1.1.} Let $\underline{X}=(X_1,X_2)\sim N_2(\theta_1,\theta_2,\sigma_1^2,\sigma_2^2,\rho)$, where $-\infty<\theta_1\leq \theta_2<\infty$, $\theta_1$ and $\theta_2$ are unknown, and $\sigma_i>0,\;i=1,2,$ and $\rho\in(-1,1)$ are known. Consider estimation of $\underline{\theta}=(\theta_1,\theta_2)$ under the sum of squared error loss functions, given by \eqref{eq:4.2}. Here $(X_1,X_2)$ is the BLEE of $(\theta_1,\theta_2)$, i.e., $c_{0,i}=0,\;i=1,2$. Also, we have $Z_i\sim N(0,\sigma_i^2),\;i=1,2,$ and $Z\sim N(0,\tau^2)$, where $\tau^2=\sigma_1^2+\sigma_2^2-2\rho \sigma_1\sigma_2$. Moreover, $f_Z(z)$ is log-concave on $\Re$. For $c_{0,1}=c_{0,2}=0$, we have
	\begin{small}
		$$\sup_{\lambda\geq 0}\alpha_1(\lambda)=\frac{p_1}{p_1+p_2}.$$
		It can be verified that
		$$ \alpha_{\infty}= \frac{p_1}{p_1+p_2}+\frac{p_2}{p_1+p_2}\lim_{\lambda\to \infty}(\alpha_1^{*}(\lambda))=\frac{p_1}{p_1+p_2}+\frac{p_2}{p_1+p_2} \lim_{\lambda\to \infty} \frac{\int_{-\infty}^{0} \lambda z\frac{1}{\tau}\phi\left(\frac{z-\lambda}{\tau}\right)dz}{\int_{-\infty}^{0}  z^2 \frac{1}{\tau}\phi\left(\frac{z-\lambda}{\tau}\right)dz}
		=-\infty .$$
	\end{small}
	
	\                Using Theorem 4.1, it follows that the estimators $\{\underline{\delta}_{\alpha}^{M}:\alpha\in (-\infty,\frac{p_1}{p_1+p_2}]\}$ are admissible within the class $\mathcal{D}=\{\underline{\delta}_{\alpha}^{M}:-\infty<\alpha<\infty\}$ (as defined by \eqref{eq:4.5} and \eqref{eq:4.6}) of isotonic regression estimators of $(\theta_1,\theta_2)$. Also, the estimators $\{\underline{\delta}_{\alpha}^{M}:\alpha>\frac{p_1}{p_1+p_2}\}$ are inadmissible and, for $\frac{p_1}{p_1+p_2}\leq \alpha_2<\alpha_1$, the estimator $\underline{\delta}_{\alpha_2}^{M}$ dominates the estimator $\underline{\delta}_{\alpha_1}^{M}$. Moreover, using Theorem 3.1, we conclude that estimators $\{\underline{\delta}_{\alpha}^{M}:\max\{0,\frac{p_1-p_2}{p_1}\}\leq \alpha<1\}$ dominate the BLEE $(X_1,X_2)$ for simultaneous estimation of $(\theta_1,\theta_2)$. For the special case $p_1=p_2=1$, the above class of admissible estimators is also obtained in Patra and Kumar (\citeyear{patra2017estimating}). Earlier Kumar and Sharma (1988) also obtained similar results for $\rho=0$ and $p_1=p_2$. For $ \alpha_0=\frac{\sigma_2(\sigma_2-\rho \sigma_1)}{\sigma_1^2+\sigma_2^2-2\rho \sigma_1 \sigma_2}$,  $p_1=\alpha_0$, $p_2=1-\alpha_0$ and $\rho<\min\{\frac{\sigma_1}{\sigma_2},\frac{\sigma_2}{\sigma_1}\}$, the restricted maximum likelihood estimator (MLE) of $\underline{\theta}$ is $\underline{\delta}_{\alpha_0}^M$ (see Patra and Kumar (2017)) and it is admissible within the class $\mathcal{D}$ of isotonic regression estimators. In this case the restricted MLE is same as the Hwang and Peddada (\citeyear{MR1272076}) estimator and Tan and Peddada (\citeyear{tan2000}) estimator.  Under $\sigma_1^2=\sigma_2^2=\sigma^2$ (say) and $p_1=p_2=1$, the restricted MLE of $\underline{\theta}$ is $\underline{\delta}_{\nu_0}^M$, where $\nu_0=\frac{1}{2}$, and it is admissible within the class $\mathcal{D}$ of isotonic regression estimators of $\underline{\theta}$. For $\rho=0$, this result is also discussed in Kumar and Sharma (1988). In general (i.e., $p_1\neq p_2$) the restricted MLE of $(\theta_1,\theta_2)$ is
	\begin{equation} \label{eq:4.11}
		\underline{\delta}_{R}(\underline{X})=\begin{cases}
			(X_1,X_2), & \text{if } X_1\leq X_2\\
			(\alpha_0X_1+(1-\alpha_0)X_2,\alpha_0 X_1+(1-\alpha_0)X_2), &\text{if } X_1>X_2
		\end{cases},
	\end{equation}
	the Hwang and Peddada (\citeyear{MR1272076}) estimator is 
	\begin{equation} \label{eq:4.12}
		\underline{\delta}_{HP}(\underline{X})=(\min\{X_1,\alpha_0 X_1+(1-\alpha_0)X_2\},\max\{X_1,\alpha_0 X_1+(1-\alpha_0)X_2\})
	\end{equation}
	and the Tan and Peddada (\citeyear{tan2000}) estimator is 
	\begin{equation} \label{eq:4.13}
		\underline{\delta}_{PDT}(\underline{X})=(\min\{X_1,\alpha_0^+ X_1+(1-\alpha_0^+)X_2\},\max\{X_1,\alpha_0^+ X_1+(1-\alpha_0^+)X_2\})
	\end{equation}
	where $\alpha_0=\frac{\sigma_2(\sigma_2-\rho \sigma_1)}{\sigma_1^2+\sigma_2^2-2\rho \sigma_1\sigma_2}$ and $\alpha_0^+=\max\{0,\alpha_0\}$. Clearly $\underline{\delta}_{R}(\underline{X}),\;\underline{\delta}_{HP}(\underline{X}),\;\underline{\delta}_{PDT}(\underline{X})\notin \mathcal{D}$, unless $\rho<\min\{\frac{\sigma_1}{\sigma_2},\frac{\sigma_2}{\sigma_1}\}$, $p_1=\alpha_0$ and $p_2=1-\alpha_0$.
	
	\vspace*{2mm}
	
	
	\vspace*{2mm}
	
	\noindent	\textbf{Example 4.1.2.} Let $X_1$ and $X_2$ be independently and identically distributed such that $X_i$ follows exponential distribution with unknown location parameter $\theta_i$ and known scale parameter $\sigma_i>0,\;i=1,2,$ where it is known apriori that $-\infty<\theta_1\leq \theta_2< \infty$. Then, $f(z_1,z_2)=f_1(z_1) f_2(z_2),\, \underline{z}=(z_1,z_2)\in\Re^2,$
	where, for known positive constants $\sigma_1>0$ and $\sigma_2>0$,
	$ f_i(z)= \frac{1}{\sigma_i} e^{-\frac{z}{\sigma_i}},\text{ if}\;z>0, i=1,2.$
	Here $c_{0,i}=\sigma_i,\;i=1,2$, and the BLEE of $(\theta_1,\theta_2)$ is $(X_1-\sigma_1,X_2-\sigma_2)$. Moreover,
	$$f_Z(z)=\begin{cases}  \frac{1}{\sigma_1+\sigma_2} e^{\frac{z}{\sigma_1}},&\text{ if} \;\; z<0  \\
		\frac{1}{\sigma_1+\sigma_2} e^{-\frac{z}{\sigma_2}},&\text{ if} \;\; z\geq 0   \end{cases}$$		
	is log-concave on $\Re$. Here, it is easy to verify that $\alpha_{\infty}=\frac{p_1}{p_1+p_2}+\frac{p_2}{p_1+p_2}\lim_{\lambda \to \infty}\alpha_1^{*}(\lambda)=-\infty$.
	Using Theorem 4.1, we conclude that the estimators $\{\underline{\delta}_{\alpha}^{M}: \alpha\in(-\infty,\frac{p_1}{p_1+p_2}] \}$ are admissible within the class $\{\underline{\delta}_{\alpha}^{M}: -\infty<\alpha<\infty \}$ of isotonic regression estimators of $(\theta_1,\theta_2)$. Moreover, for $\frac{p_1}{p_1+p_2}\leq \alpha_2<\alpha_1<\infty$, the isotonic regression estimator $ \underline{\delta}_{\alpha_2}^{M}$ dominates the estimator $\underline{\delta}_{\alpha_1}^{M}$. In particular the BLEE $\underline{\delta}_{0}(\underline{X})=(X_1-\sigma_1,X_2-\sigma_2)$ is inadmissible for estimating $(\theta_1,\theta_2)$ and is dominated by the isotonic regression estimators $\underline{\delta}_{\alpha}^{M}(\underline{X})=(\delta_{1,\alpha}^{M}(\underline{X}),\delta_{2,\alpha}^{M}(\underline{X})),\; \frac{p_1}{p_1+p_2}\leq \alpha<1$. Also, using Theorem 3.1 the isotonic regression estimators $\underline{\delta}_{\alpha}^M$, for $\max\{0,\frac{p_1-p_2}{p_1}\}\leq \alpha<  1$, dominate the BLEE $(X_1-\sigma_1,X_2-\sigma_2)$.
	

	\subsection{\textbf{Simulation Study For Estimation of Location Parameter $(\theta_1,\theta_2)$}}
	\label{sec:3.3}
	\noindent

	\vspace*{3mm}

	In Example 4.1.1, we considered estimation of the location parameters $(\theta_1,\theta_2)$ of a bivariate normal distribution having order restricted location parameters (i.e., $\theta_1\leq \theta_2$), known variances and known correlation coefficient ($\sigma_1>0$, $\sigma_2>0$ and $\rho\in (-1,1)$). We have shown that the class of isotonic regression estimators (IREs) $\{\underline{\delta}_{\alpha}^{M}: \alpha\in(-\infty,\frac{p_1}{p_1+p_2}] \}$ is admissible within the class $\mathcal{D}=\{\underline{\delta}_{\alpha}^{M}: -\infty<\alpha<\infty \}$ and also, estimators in the class $\{\underline{\delta}_{\alpha}^{M}: \max\{0,\frac{p_1-p_2}{p_1}\}\leq \alpha<1 \}$ dominate the BLEE $\underline{\delta}_{0}(\underline{X})=(X_1,X_2)$. To further evaluate the performances of various estimators under the loss function $	L(\underline{\theta},\underline{a})=\frac{1}{\sigma_1^2} (a_1-\theta_1)^2 +\frac{1}{\sigma_2^2} (a_2-\theta_2)^2,\;\underline{\theta}=(\theta_1,\theta_2)\in\Theta_0,\;\underline{a}=(a_1,a_2)\in\Re^2$ (i.e., $p_1=\frac{1}{\sigma_1^2}$ and $p_2=\frac{1}{\sigma_2^2}$), in this section, we compare the risk performances of the BLEE $(X_1,X_2)$, the isotonic regression estimators (IREs) $\underline{\delta}_{\alpha}^{M}(\underline{X})=(\delta_{1,\alpha}^{M}(\underline{X}),\delta_{2,\alpha}^{M}(\underline{X}))$ with $\alpha=\frac{p_1}{p_1+p_2}=\frac{\sigma_2^2}{\sigma_1^2+\sigma_2^2}$, the restricted MLE $\underline{\delta}_{R}(\underline{X})$ (as defined by \eqref{eq:4.11}), the Hwang and Peddada (\citeyear{MR1272076}) estimator $\underline{\delta}_{HP}(\underline{X})$ (as defined by \eqref{eq:4.12}) and the Tan and Peddada (\citeyear{tan2000}) estimator $\underline{\delta}_{PDT}(\underline{X})$ (as defined by \eqref{eq:4.13}), numerically, through the Monte Carlo simulations. Note that, for $\rho<\min\{\frac{\sigma_1}{\sigma_2},\frac{\sigma_2}{\sigma_1}\}$, the restricted MLE $\underline{\delta}_{R}(\underline{X})$, the Hwang and Peddada (\citeyear{MR1272076}) estimator $\underline{\delta}_{HP}(\underline{X})$ and the Tan and Peddada (\citeyear{tan2000}) estimator $\underline{\delta}_{PDT}(\underline{X})$ are same.
	
	\vspace*{2mm}

	For simulations, we generated 50000 samples of size 1 each from relevant bivariate distributions and computed the simulated risks of various estimators.
	

	\FloatBarrier
	\begin{figure}[h!]
		\begin{subfigure}{0.48\textwidth}
			\includegraphics[width=80mm,scale=1.2]{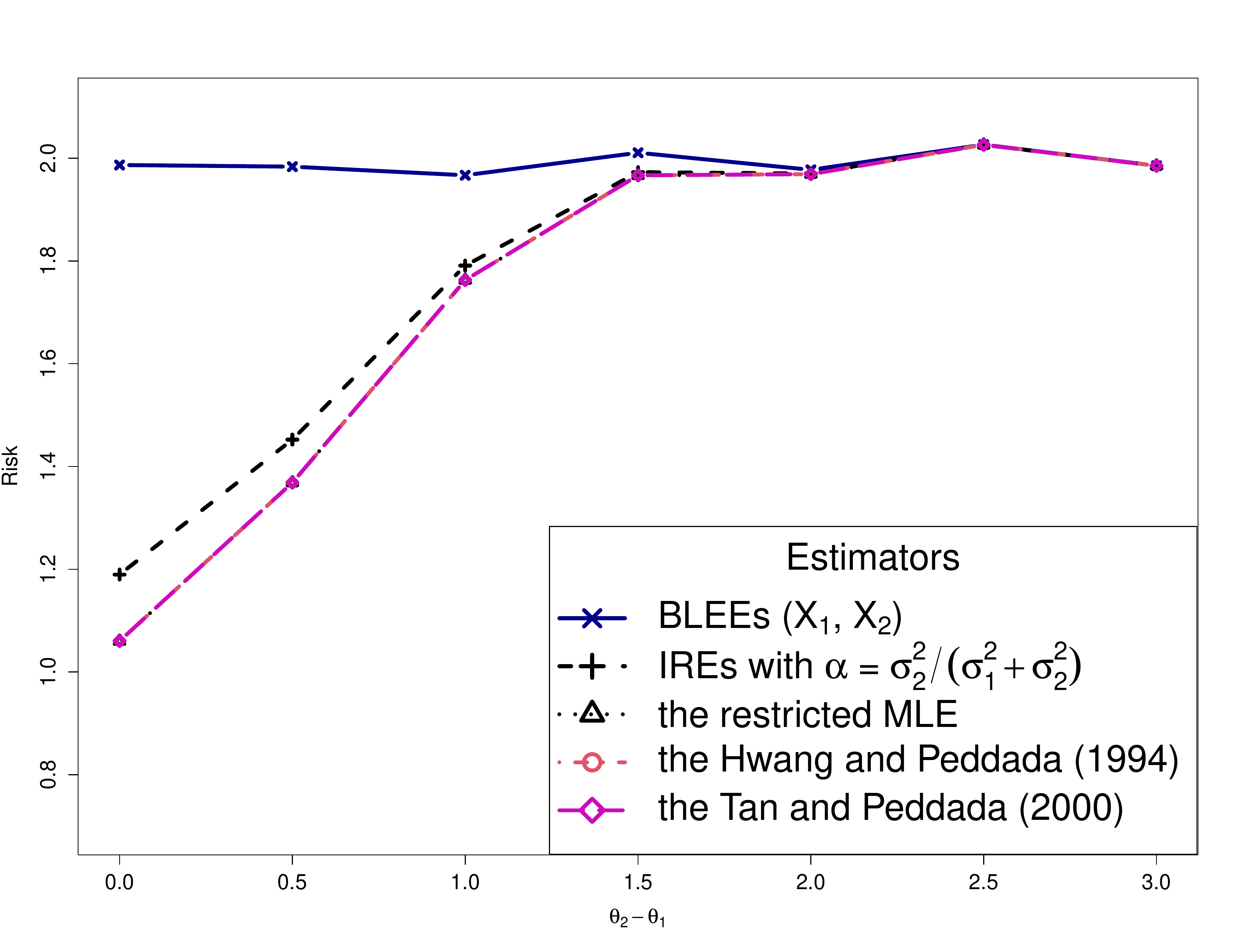} 
			\caption{$\sigma_1=0.2$, $\sigma_2=0.5$ and $\rho=-0.9$.} 
			\label{fig7:a} 
		\end{subfigure}
		\begin{subfigure}{0.48\textwidth}
			\includegraphics[width=80mm,scale=1.2]{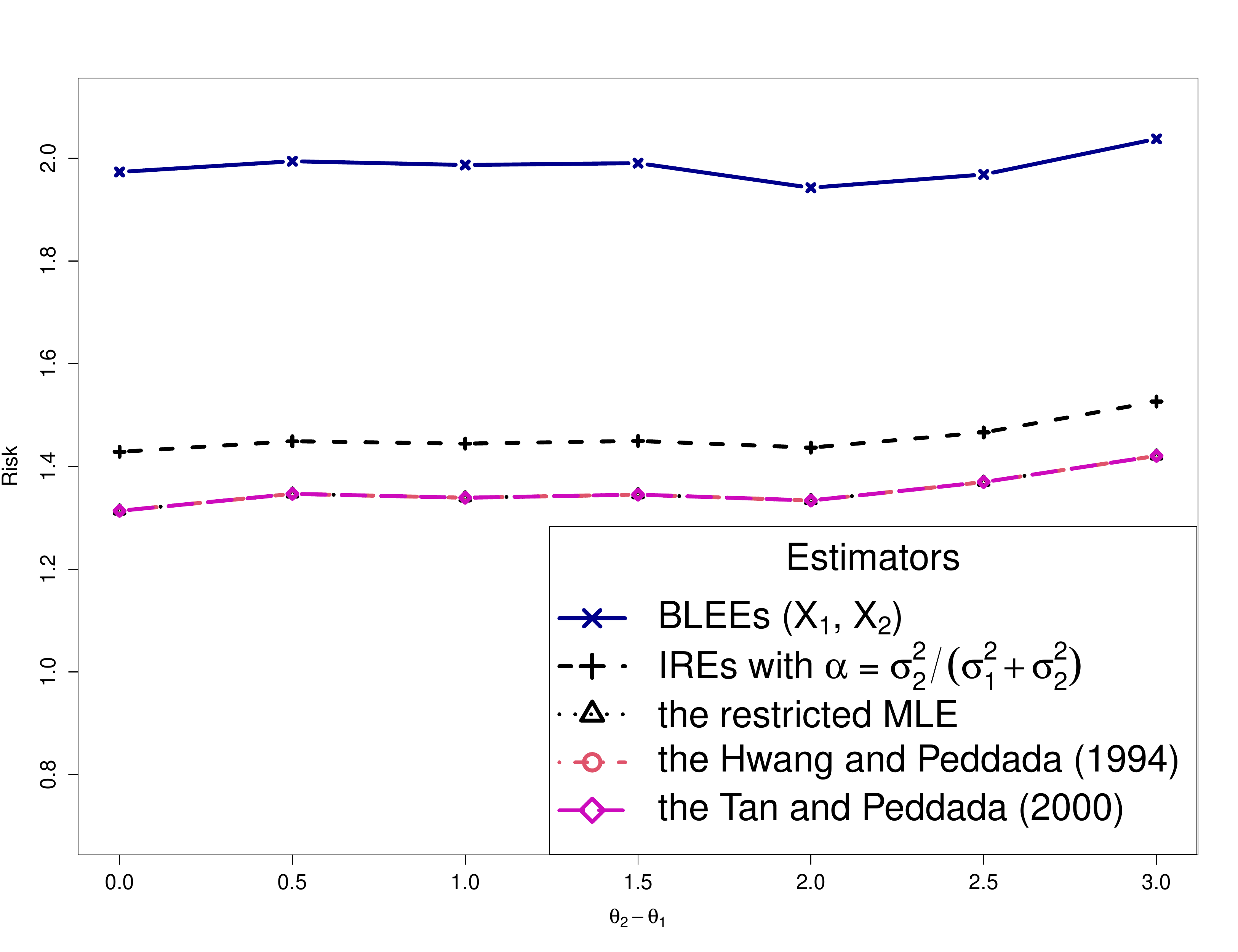} 
			\caption{$\sigma_1=10$, $\sigma_2=1$ and $\rho=-0.5$.} 
			\label{fig7:b} 
		\end{subfigure}
		\\	\begin{subfigure}{0.48\textwidth}
			
			\includegraphics[width=80mm,scale=1.2]{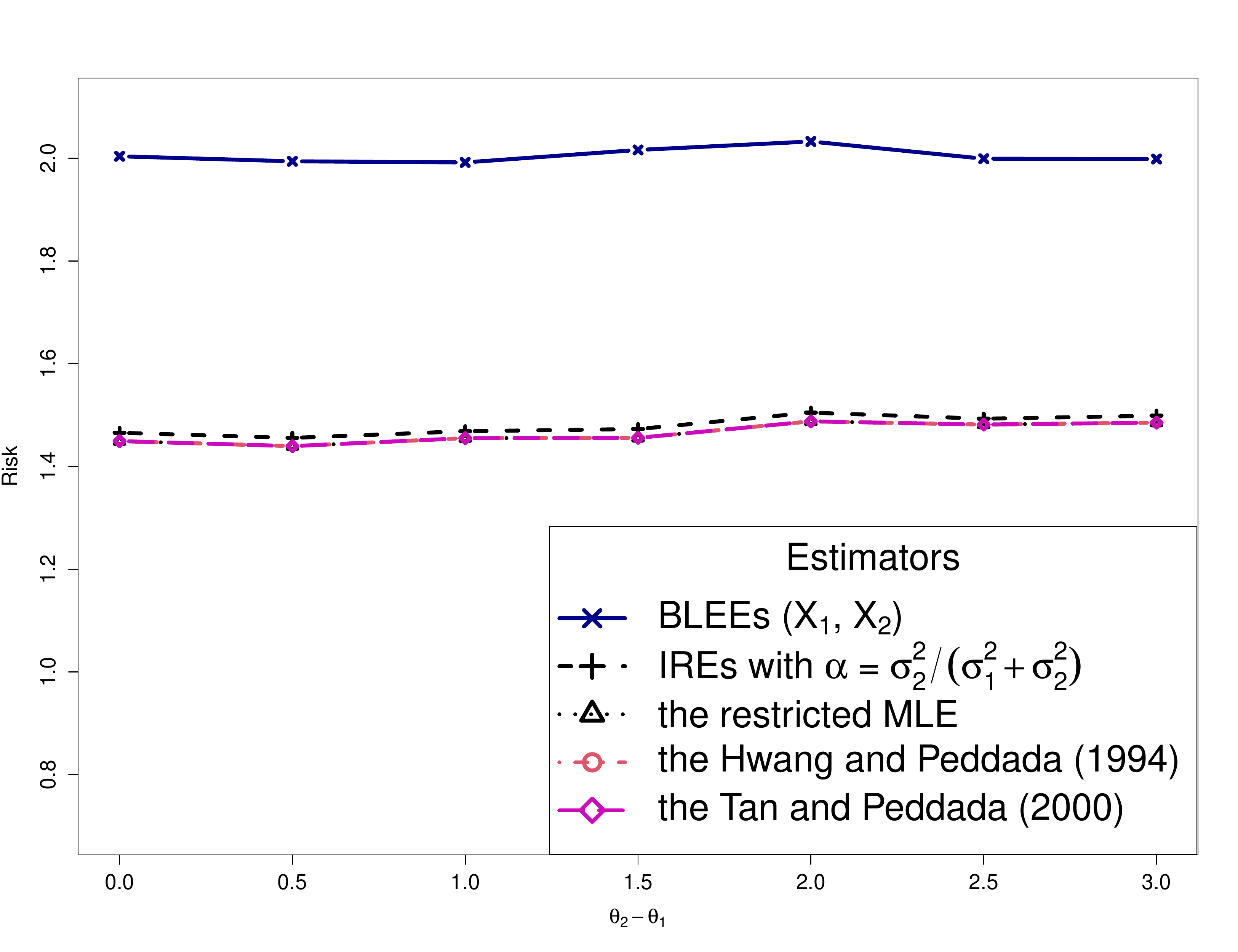} 
			
			\caption{$\sigma_1=2$, $\sigma_2=10$ and $\rho=-0.2$.} 
			\label{fig7:c} 
		\end{subfigure}
		\begin{subfigure}{0.48\textwidth}
			\includegraphics[width=80mm,scale=1.2]{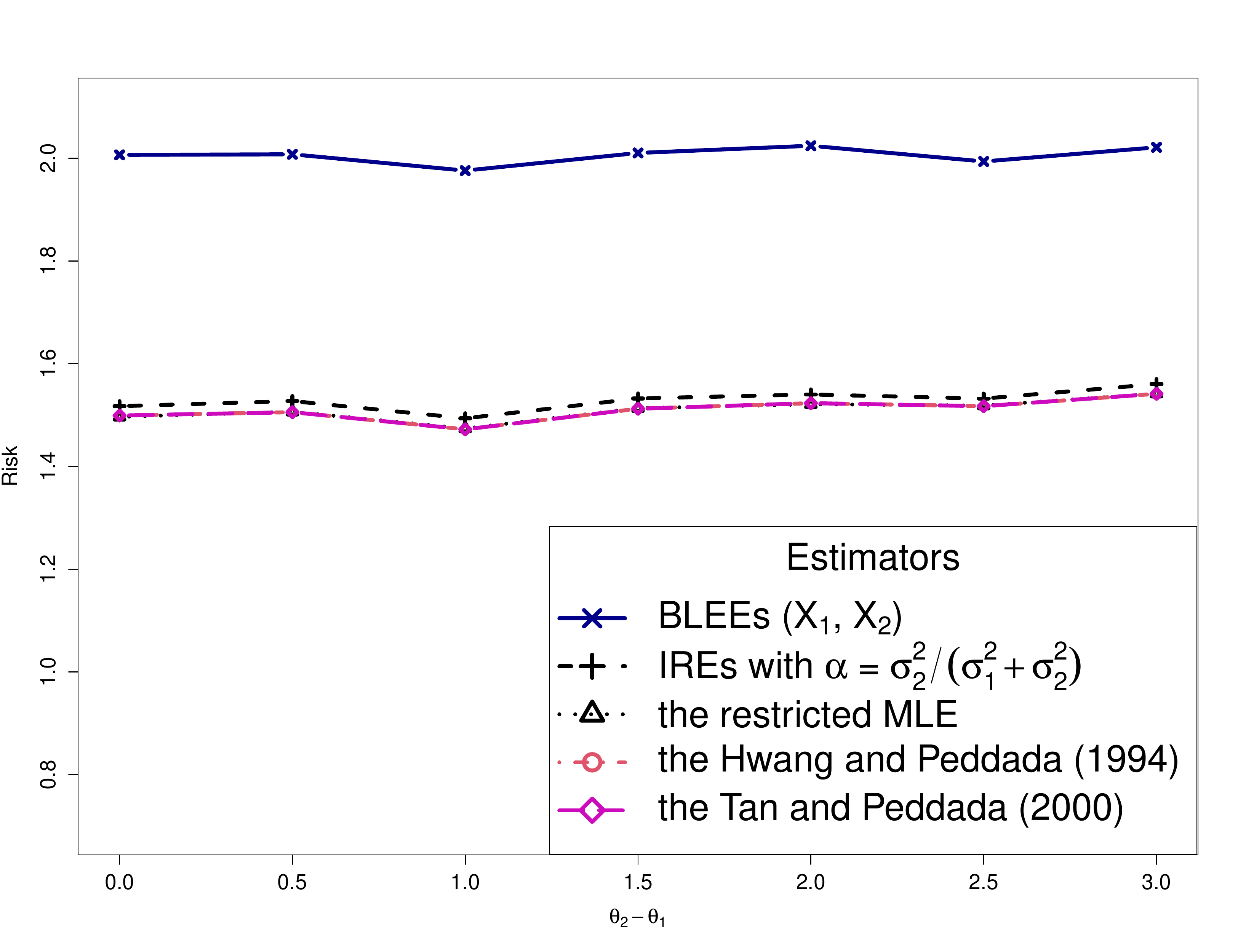} 
			\caption{ $\sigma_1=1$, $\sigma_2=10$ and $\rho=0.2$.} 
			\label{fig7:d}  
		\end{subfigure}
		\\	\begin{subfigure}{0.48\textwidth}
			\centering
			
			\includegraphics[width=80mm,scale=1.2]{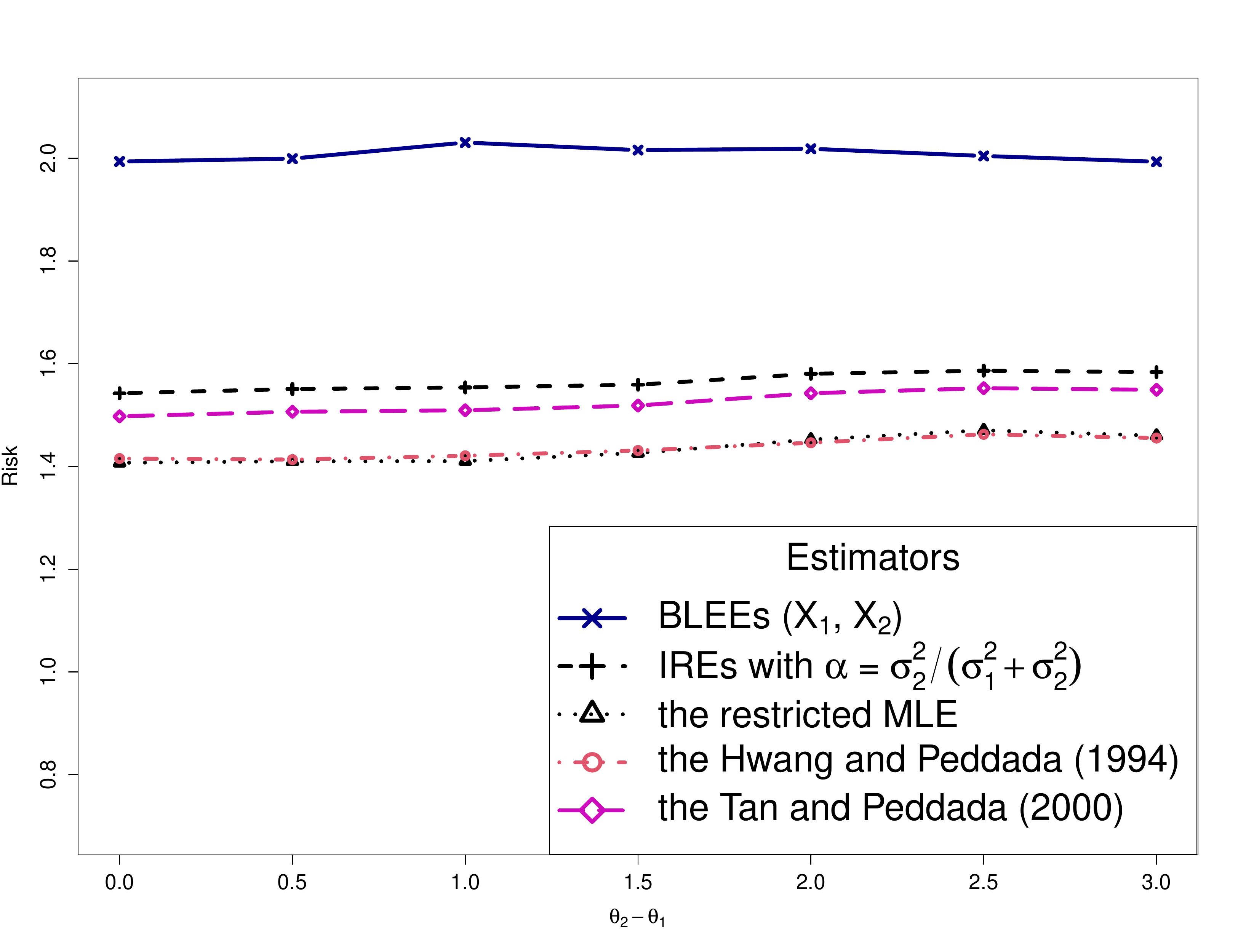} 
			\caption{$\sigma_1=10$, $\sigma_2=1$ and $\rho=0.5$.} 
			\label{fig7:e} 
		\end{subfigure}
		\begin{subfigure}{0.48\textwidth}
			\centering
			\includegraphics[width=80mm,scale=1.2]{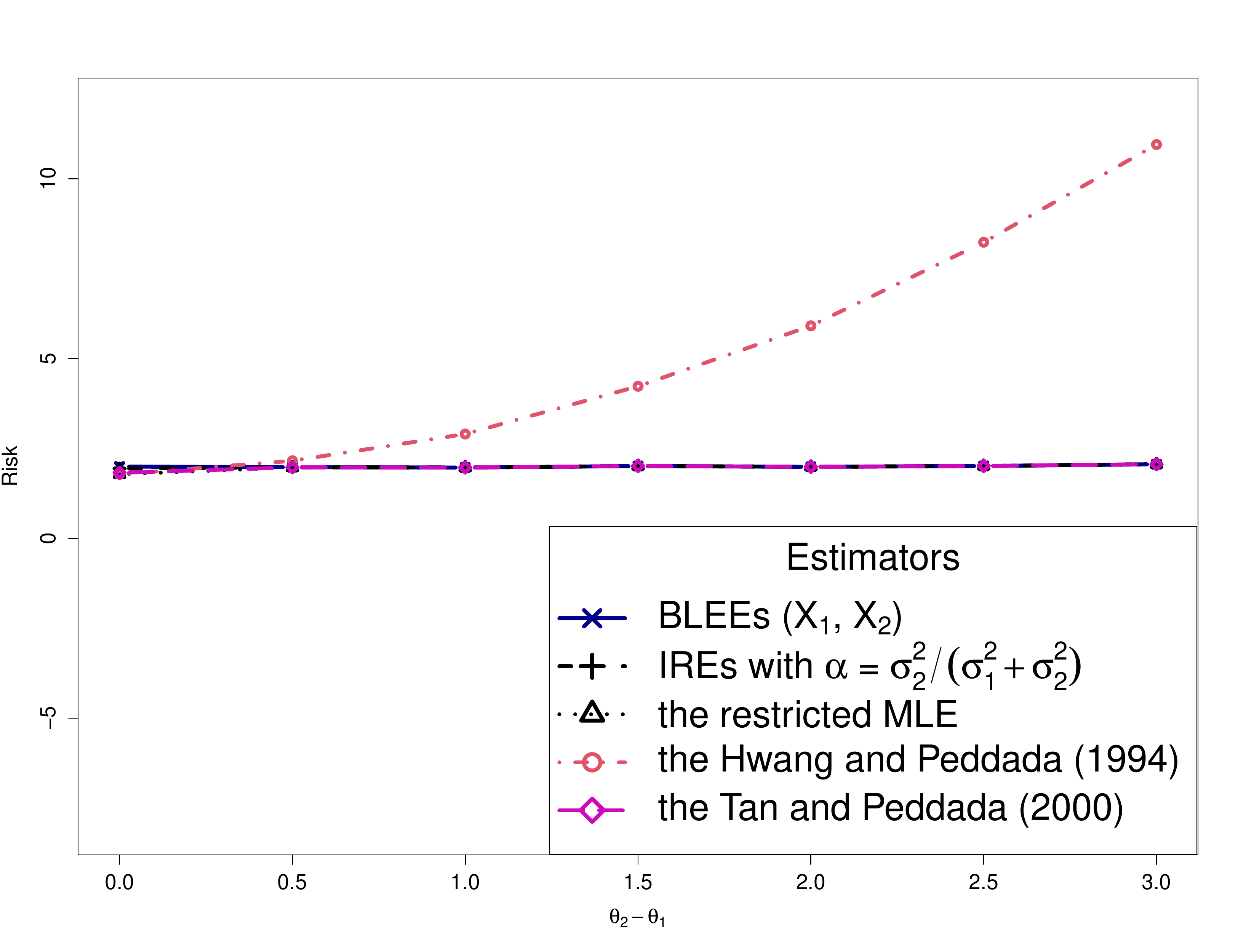} 
			\caption{ $\sigma_1=0.5$, $\sigma_2=0.4$ and $\rho=0.9$.} 
			\label{fig7:f}  
		\end{subfigure}
		
		\caption{Risk plots of the BLEE $(X_1,X_2)$, the IREs with $\alpha=\frac{\sigma_2^2}{\sigma_1^2+\sigma_2^2}$ and the restricted MLE.}
		\label{fig1}
	\end{figure}
	\FloatBarrier
	
	\noindent

	\vspace*{2mm}
	
	The simulated values of risks of various estimators are plotted in Figure 1. The following observations are evident from Figure 1:

	\noindent (i) The IREs with $\alpha=\frac{\sigma_2^2}{\sigma_1^2+\sigma_2^2}$ and the restricted MLE always dominate the BLEE $(X_1,X_2)$.\vspace*{2mm}
	
	\noindent (ii) The Restricted MLE always perform better than the other estimators.

	\vspace*{7mm}

	Under the sum of squared error loss functions \eqref{eq:4.2}, in Example 4.1.2, we considered estimation of the location parameters $(\theta_1,\theta_2)$ of two independent exponential distributions having unknown order restricted location parameters (i.e., $\theta_1\leq \theta_2$) and known scale parameters ($\sigma_1>0$ and $\sigma_2>0$). We have shown that the class of isotonic regression estimators (IREs) $\{\underline{\delta}_{\alpha}^{M}: \alpha\in(-\infty,\frac{p_1}{p_1+p_2}] \}$ is admissible within the class $\mathcal{D}=\{\underline{\delta}_{\alpha}^{M}: -\infty<\alpha<\infty \}$. Also, the class of estimators $\{\underline{\delta}_{\alpha}^{M}: \max\{0,\frac{p_1-p_2}{p_1}\}\leq \alpha<1 \}$ dominates the BLEE $\underline{\delta}_{0}(\underline{X})=(X_1-\sigma_1,X_2-\sigma_2)$. To further evaluate the performances of various estimators under the loss function \eqref{eq:4.2} with $p_1=p_2=1$, in this section, we compare the risk performances of the isotonic regression estimators (IREs) $\underline{\delta}_{1}^{M}(\underline{X})=(\delta_{1,1}^{M}(\underline{X}),\delta_{2,1}^{M}(\underline{X}))=(X_1-\sigma_1,X_2-\sigma_2)$, $\underline{\delta}_{0.75}^{M}(\underline{X})=(\delta_{1,0.75}^{M}(\underline{X}),\delta_{2,0.75}^{M}(\underline{X}))$, $\underline{\delta}_{0.5}^{M}(\underline{X})=(\delta_{1,0.5}^{M}(\underline{X}),\delta_{2,0.5}^{M}(\underline{X}))$, $\underline{\delta}_{0}^{M}(\underline{X})=(\delta_{1,0}^{M}(\underline{X}),\delta_{2,0}^{M}(\underline{X}))$ and the restricted MLE $\underline{\delta}_R(\underline{X})=(\min\{X_1,X_2\},X_2)$, numerically, through the Monte Carlo simulations; here 
	$$\delta_{1,\alpha}^{M}(\underline{X})=\begin{cases}
		X_1-\sigma_1,&\text{ if} \;\; X_1-\sigma_1\leq X_2-\sigma_2\\
		\alpha(X_1-\sigma_1) + (1-\alpha)(X_2-\sigma_2), &\text{ if} \; \; X_1-\sigma_1> X_2-\sigma_2
	\end{cases},$$
	$$\text{and   } \delta_{2,\alpha}^{M}(\underline{X})=\begin{cases}
		X_2-\sigma_2,&\text{ if} \;\; X_1-\sigma_1\leq X_2-\sigma_2\\
		(1-\alpha)(X_1-\sigma_1) + \alpha(X_2-\sigma_2), &\text{ if} \; \; X_1-\sigma_1> X_2-\sigma_2
	\end{cases}.$$

	For simulations, we generated 50000 samples of size 1 each from relevant exponential distributions and computed the simulated risks of estimators $\underline{\delta}_{1}^{M}(\underline{X})=(X_1-\sigma_1,X_2-\sigma_2)$, $\underline{\delta}_{0.75}^{M}(\underline{X})$, $\underline{\delta}_{0.5}^{M}(\underline{X})$, $\underline{\delta}_{0}(\underline{X})$ and $\underline{\delta}_{R}(\underline{X})$.
	\vspace*{2mm}
	
	The simulated values of risks of various estimators are plotted in Figure 2. The following observations are evident from Figure 2:
	\\~\\(i) The risk function values of the IREs $\underline{\delta}_{0.75}^{M}(\underline{X})$, $\underline{\delta}_{0.5}^{M}(\underline{X})$, $\underline{\delta}_{0.25}^{M}(\underline{X})$ and $\underline{\delta}_{0}^{M}(\underline{X})$ are nowhere larger than the risk function values of the BLEE $(X_1-\sigma_1,X_2-\sigma_2)$, which is in 
	conformity with theoretical findings of Example 4.1.2.
	\noindent \\~\\(ii) From Figure 2, we can observe that the IREs $\underline{\delta}_{0.5}^{M}(\underline{X})$, $\underline{\delta}_{0}^{M}(\underline{X})$ and $\underline{\delta}_{R}^{M}(\underline{X})$ are not comparable and the IREs $\underline{\delta}_{1}^{M}(\underline{X})$ and $\underline{\delta}_{0.75}^{M}(\underline{X})$ are inadmissible. This is in conformity with theoretical findings of Example 4.1.2.
	\\~\\(iii) When $\sigma_1>>\sigma_2$, for small and moderate values of $\theta_2-\theta_1$, the restricted MLE outperforms the other estimators. 
	\\~\\(iv) There is no clear cut winner between various estimators, but performance of the BLEE $\underline{\delta}_{1}^{M}(\underline{X})=(X_1-\sigma_1,X_2-\sigma_2)$ and $\underline{\delta}_{0.75}^{M}(\underline{X})$ are worse than other estimators. Also, when $\sigma_1\leq \sigma_2$, the performance of the restricted MLE is the worst among other estimators.

	\FloatBarrier
	\begin{figure}[h!]
		\begin{subfigure}{0.48\textwidth}
			\centering
			\includegraphics[width=85mm,scale=1.2]{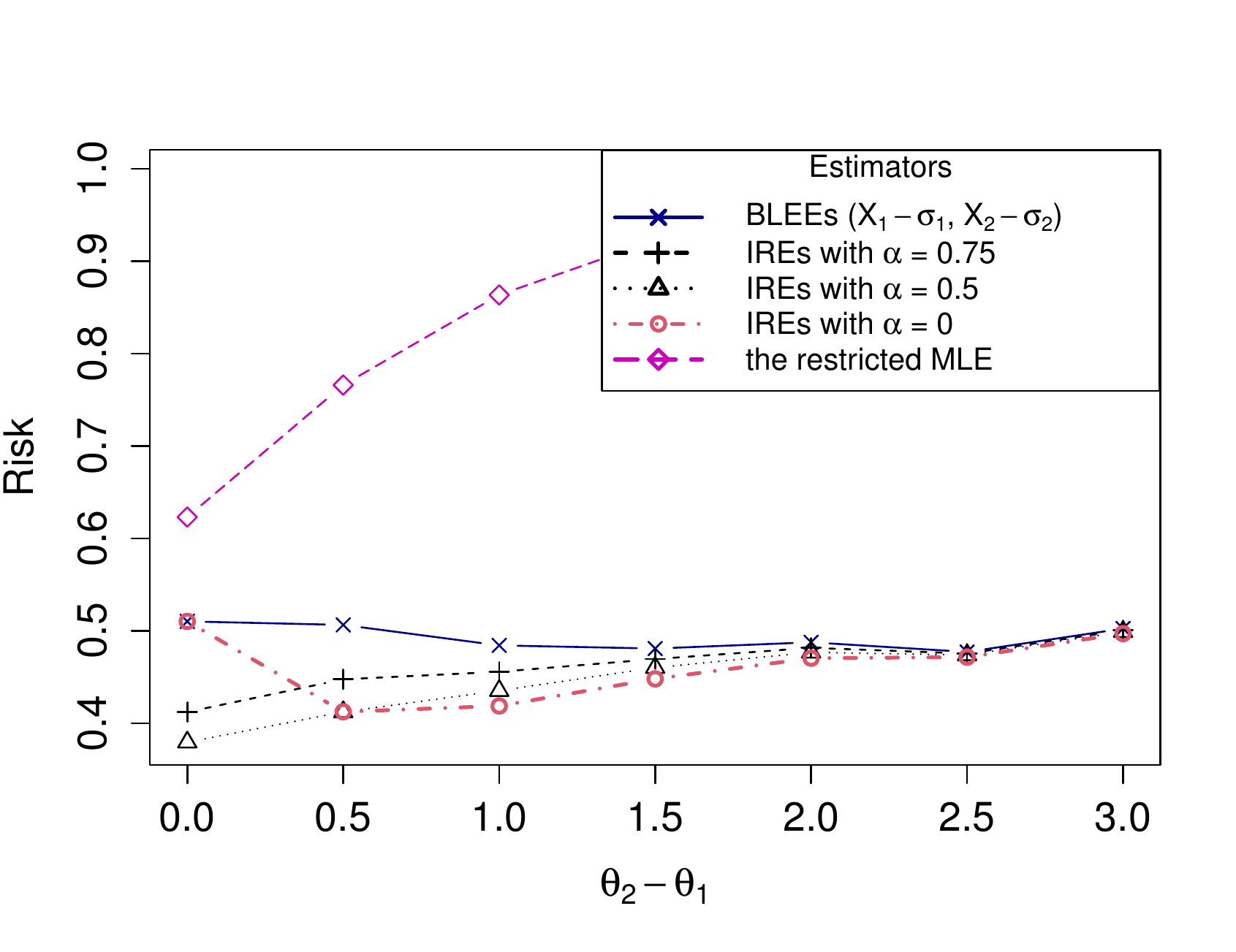} 
			\caption{$\sigma_1=0.5$ and $\sigma_2=0.5$.} 
			\label{fig1:a} 
		\end{subfigure}
		\begin{subfigure}{0.48\textwidth}
			\centering
			\includegraphics[width=85mm,scale=1.2]{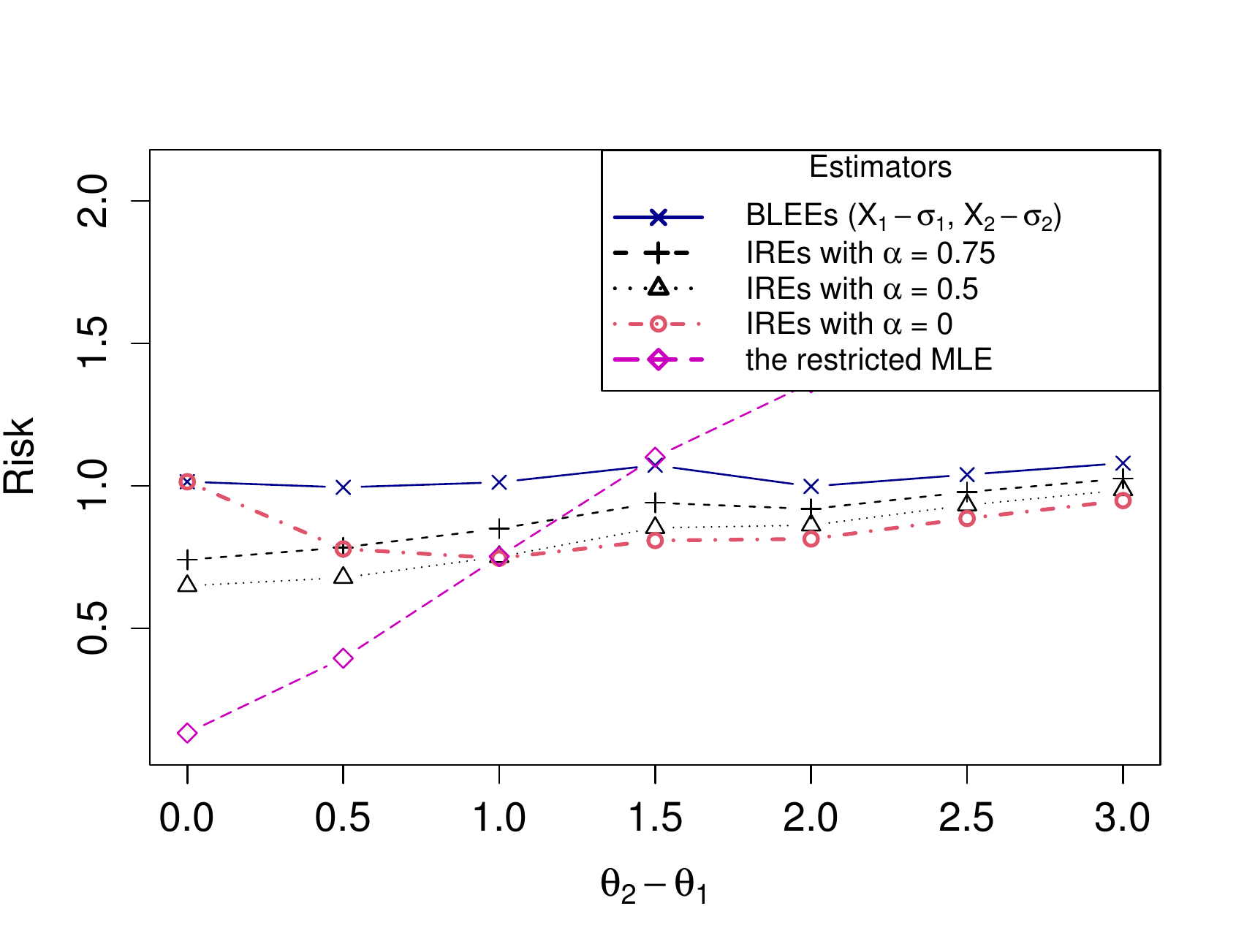} 
			\caption{$\sigma_1=1$ and $\sigma_2=0.2$.} 
			\label{fig1:b} 
		\end{subfigure}
		\\	\begin{subfigure}{0.48\textwidth}
			\centering
			
			\includegraphics[width=85mm,scale=1.2]{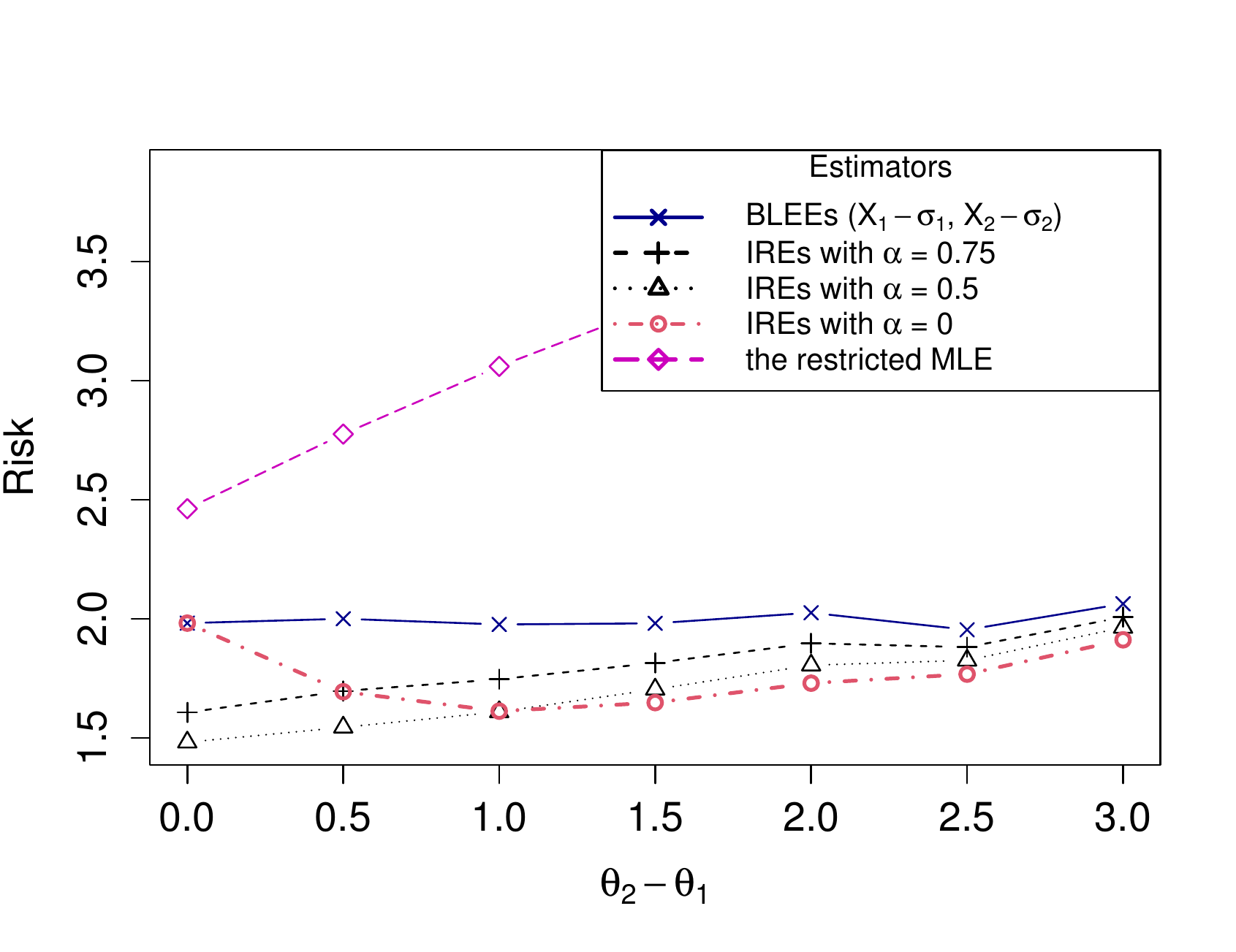} 
			
			\caption{$\sigma_1=1$ and $\sigma_2=1$.} 
			\label{fig1:c} 
		\end{subfigure}
		\begin{subfigure}{0.48\textwidth}
			\centering
			\includegraphics[width=85mm,scale=1.2]{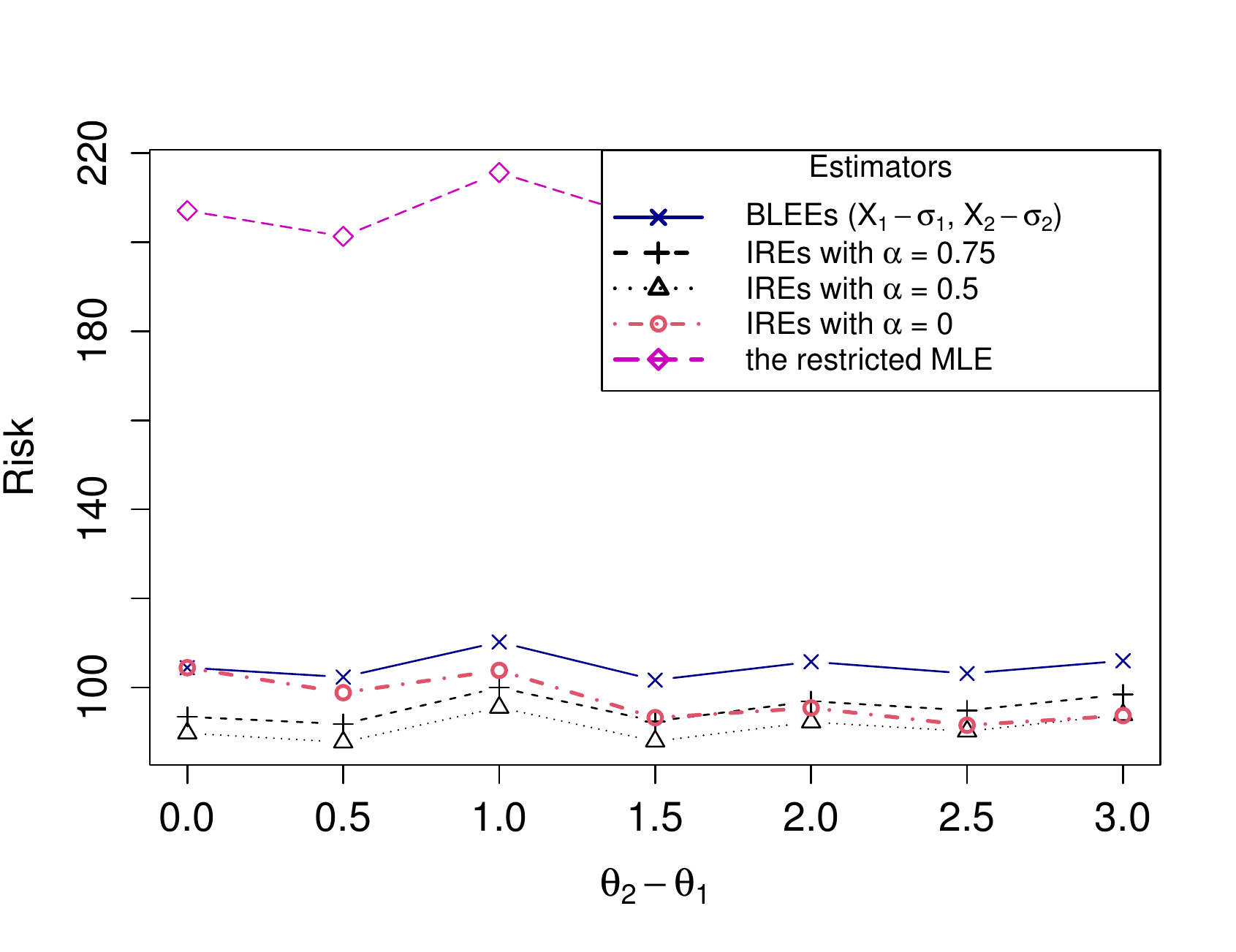} 
			\caption{ $\sigma_1=2$ and $\sigma_2=10$.} 
			\label{fig1:d}  
		\end{subfigure}
		\\	\begin{subfigure}{0.48\textwidth}
			\centering
			
			\includegraphics[width=85mm,scale=1.2]{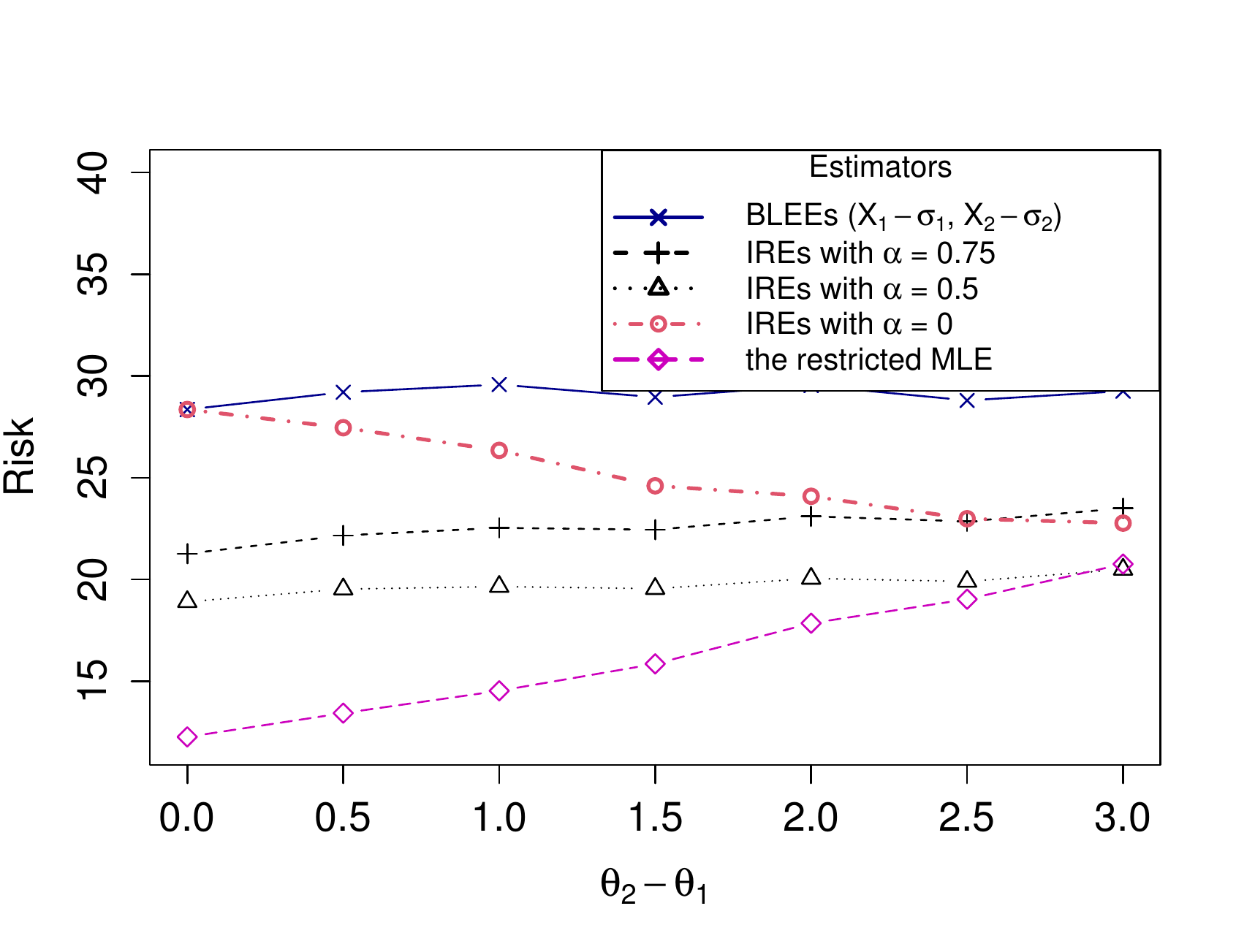} 
			
			\caption{$\sigma_1=5$ and $\sigma_2=2$.} 
			\label{fig1:e} 
		\end{subfigure}
		\begin{subfigure}{0.48\textwidth}
			\centering
			\includegraphics[width=85mm,scale=1.2]{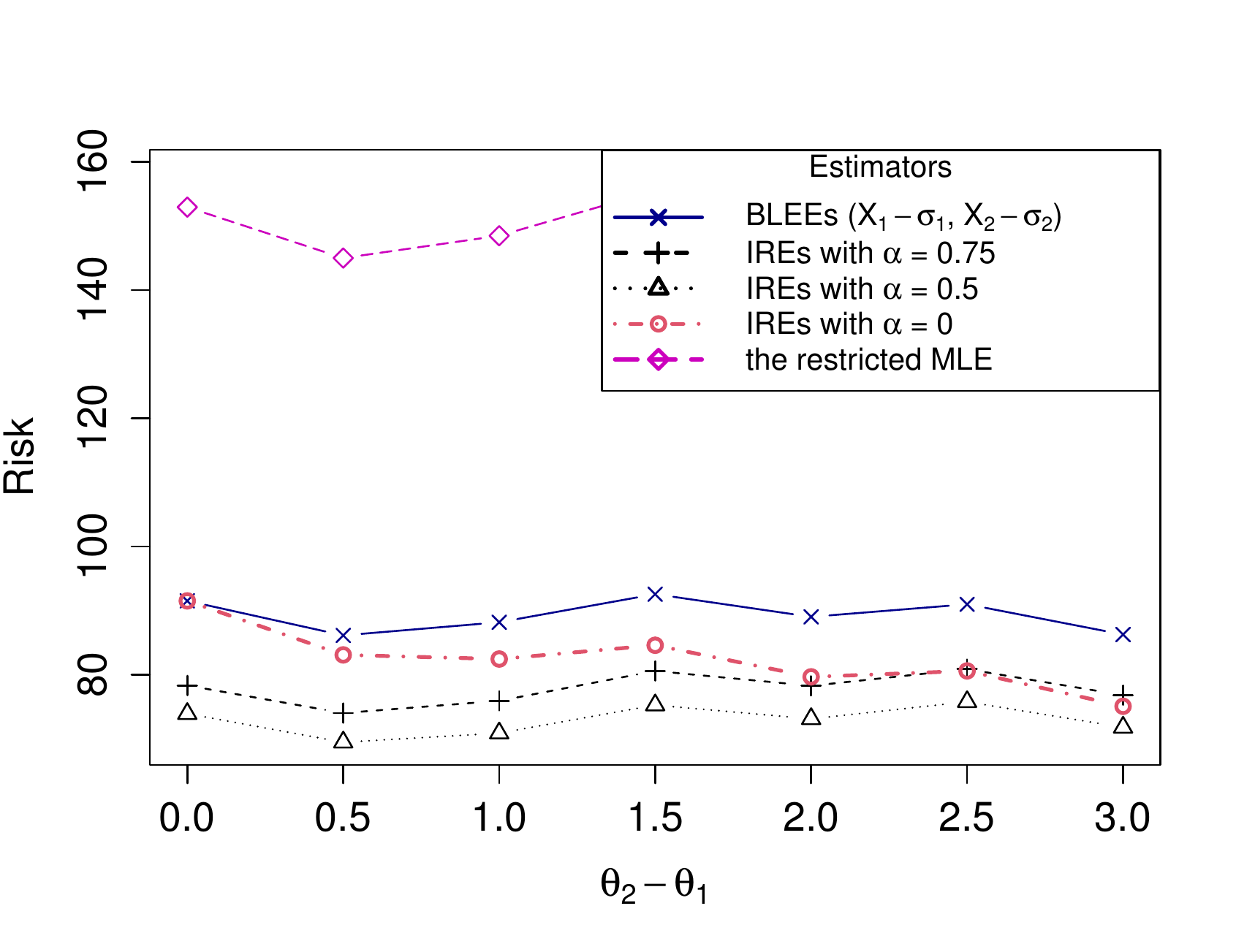} 
			\caption{ $\sigma_1=5$ and $\sigma_2=8$.} 
			\label{fig1:f}  
		\end{subfigure}
		\caption{Risk plots of estimators $\underline{\delta}_{1}^{M}(\underline{X})=(X_1-\sigma_1,X_2-\sigma_2)$, $\underline{\delta}_{0.75}^{M}(\underline{X})$, $\underline{\delta}_{0.5}^{M}(\underline{X})$, $\underline{\delta}_{0}^{M}(\underline{X})$ and $\underline{\delta}_{R}^{M}(\underline{X})$ against the values of $\lambda=\theta_2-\theta_1$.}
		\label{fig2}
	\end{figure}
	\FloatBarrier

	\section{\textbf{Admissibility of Isotonic regression estimators of Scale Parameters under weighted sum of squared errors loss}}

	Let the random vector $\underline{X}=(X_1,X_2)$ have the Lebesgue p.d.f. belonging to the scale family
	\begin{equation}\label{eq:5.1}
		f_{\underline{\theta}}(x_1,x_2)=	\frac{1}{\theta_1 \theta_2}f\!\left(\frac{x_1}{\theta_1},\frac{x_2}{\theta_2}\right),\; (x_1,x_2)\in \Re^2, \;\underline{\theta}=(\theta_1,\theta_2)\in \Theta_0,
	\end{equation} 
	where $\Theta_0=\{(\theta_1,\theta_2)\in \Re_{++}^2: \theta_1\leq \theta_2\}$. Throughout, we will assume that the distributional support of $\underline{X}=(X_1,X_2)$ is a subset of $\Re_{++}^2$. Generally, $\underline{X}=(X_1,X_2)$ would be a minimal-sufficient statistic based on a bivariate random sample or two independent random samples.
	\vspace*{2mm}
	
	\noindent For estimating $\underline{\theta}=(\theta_1,\theta_2)$, consider the loss function
	\begin{equation}\label{eq:5.2}
		L(\underline{\theta},\underline{a})=p_1 (a_1-\theta_1)^2 +p_2 (a_2-\theta_2)^2, \;\underline{\theta}=(\theta_1,\theta_2)\in\Theta_0,\;\underline{a}=(a_1,a_2)\in\mathcal{A}=\Re_{++}^2,
	\end{equation}
	where $p_1>0$ and $p_2>0$ are pre-specified constants.
	\vspace*{2mm}

	Under the (unrestricted) parameter space $\Theta=\Re_{++}^2$, the problem of simultaneously estimating $\theta_1$ and $\theta_2$, with the loss function \eqref{eq:5.2}, is invariant under the multiplicative group of transformations $\mathcal{G}=\{g_{c_1,c_2}:(c_1,c_2)\in\Re_{++}^2 \}$, where $g_{c_1,c_2}(x_1,x_2)=(c_1x_1,c_2x_2)$, $(x_1,x_2)\in\Re_{++}^2,\;  (c_1,c_2)\in\Re_{++}^2$. The group $\mathcal{G}$ induces the groups of transformations $\overline{\mathcal{G}}=\{\overline{g}_{c_1,c_2}:(c_1,c_2)\in\Re_{++}^2 \}$ and $\tilde{\mathcal{G}}=\{\tilde{g}_{c_1,c_2}:(c_1,c_2)\in\Re_{++}^2 \}$ on the parameter space $\Theta$ and the action space $\mathcal{A}$, respectively, where $\overline{g}_{c_1,c_2}(\theta_1,\theta_2)=(c_1\theta_1,c_2 \theta_2)$, $\tilde{g}_{c_1,c_2}(a_1,a_2)=(c_1a_1,c_2a_2)$, $(\theta_1,\theta_2)\in\Theta,\;(a_1,a_2)\in\mathcal{A}=\Re_{++}^2,\;  (c_1,c_2)\in\Re_{++}^2$. Under the group of transformations $\mathcal{G}$, an estimator $\underline{\delta}=(\delta_1,\delta_2)$ is invariant (scale invariant) for estimating $\underline{\theta}=(\theta_1,\theta_2)$ if, and only if
	\begin{equation*}
		\underline{\delta}=(\delta_1(\underline{X}),\delta_2(\underline{X}))=(c_1X_1,c_2X_2),
	\end{equation*}
	for some $\underline{c}=(c_1,c_2)\in \Re_{++}^2$. The best scale equivariant estimator (BSEE) of $(\theta_1,\theta_2)$ is 
	\begin{equation*}
		\underline{\delta}_{0}(\underline{X})=(\delta_{0,1}(\underline{X}),\delta_{0.2}(\underline{X}))=(c_{0,1}X_1,c_{0,2}X_2),
	\end{equation*}
	where, for $Z_i=\frac{X_i}{\theta_i},$ $c_{0,i}=\frac{E[Z_i]}{E[Z_i^2]},\;i=1,2.$

	\vspace*{3mm}

	\noindent Based on BSEEs $(c_{0,1}X_1,c_{0,2}X_2)$, define isotonic regression estimators $\underline{\delta}^M_{\alpha}(\underline{X})=(\delta_{1,\alpha}^{M}(\underline{X}),\delta_{2,\alpha}^{M}(\underline{X}))$, where
	\begin{align} \label{eq:5.3}
		\delta_{1,\alpha}^{M}(\underline{X})=\begin{cases}
			c_{0,1}X_1, &\text{ if}\quad c_{0,1}X_1\leq c_{0,2}X_2 \\
			\alpha c_{0,1}X_1+(1-\alpha)c_{0,2}X_2,&\text{ if}\quad c_{0,1}X_1> c_{0,2}X_2
		\end{cases} 
	\end{align}
	and       \begin{align}\label{eq:5.4}
		\delta_{2,\alpha}^{M}(\underline{X})=\begin{cases}
			c_{0,2}X_2, &\text{ if}\quad c_{0,1}X_1\leq c_{0,2}X_2 \\
			\frac{p_1}{p_2}(1-\alpha) c_{0,1}X_1+(1-\frac{p_1}{p_2}(1-\alpha))c_{0,2}X_2,&\text{ if}\quad c_{0,1}X_1> c_{0,2} X_2
		\end{cases},\;\;\alpha\in \Re.
	\end{align}

	Let $\mathcal{D}=\{\underline{\delta}_{\alpha}^{M}=(\delta_{1,\alpha}^{M}(\underline{X}),\delta_{2,\alpha}^{M}(\underline{X})):\alpha\in\Re\}$ be the class of isotonic regression estimators (IREs) based on BSEEs $(c_{0,1}X_1,c_{0,2}X_2)$. The risk function of the estimator $\underline{\delta}_{\alpha}^{M}$, $\alpha\in \Re$, is given by
	
	\begin{small}
		\begin{align*}
			R(\underline{\theta},\underline{\delta}_{\alpha}^{M}) \nonumber
			&=p_1E_{\underline{\theta}}\left[(c_{0,1}\theta_1 Z_1-\theta_1)^2\, I_{\left(\frac{c_{0,1}}{c_{0,2}\lambda},\infty\right)}(Z)\right]\\ \nonumber
			&\quad +p_1 E_{\underline{\theta}}\left[\{\alpha (c_{0,1}\theta_1 Z_1-c_{0,2}\theta_2 Z_2)+c_{0,2}\theta_2 Z_2-\theta_1\}^2 \, I_{\left(0,\frac{c_{0,1}}{c_{0,2}\lambda}\right)}(Z)\right]\\ \nonumber
			&\quad +p_2 E_{\underline{\theta}}\left[(c_{0,2}\theta_2 Z_2-\theta_2)^2\, I_{\left(\frac{c_{0,1}}{c_{0,2}\lambda},\infty\right)}(Z)\right]\\
			&\quad + p_2 E_{\underline{\theta}}\left[\bigg\{\alpha \frac{p_1}{p_2}(c_{0,2}\theta_2 Z_2-c_{0,1}\theta_1 Z_1)-\frac{p_1-p_2}{p_2}c_{0,2}\theta_2 Z_2+\frac{p_1}{p_2}c_{0,1}\theta_1 Z_1-\theta_2\bigg\}^2 \, I_{\left(0,\frac{c_{0,1}}{c_{0,2}\lambda}\right)}(Z)\right],\; \underline{\theta}\in \Theta_0,
		\end{align*}
	\end{small}
	\noindent where $Z_1=\frac{X_1}{\theta_1}$, $Z_2=\frac{X_2}{\theta_2}$ and $Z=\frac{Z_2}{Z_1}$. Define $\lambda=\frac{\theta_2}{\theta_1}$ ($\lambda\geq 1$) and let $f_Z(\cdot)$ denote the p.d.f. of $Z=\frac{Z_2}{Z_1}$, so that $f_Z(z)=\int_{-\infty}^{\infty} \vert s \vert f(s,sz)\,ds,\; -\infty<z<\infty$. 
	Clearly, for any fixed $\underline{\theta}\in\Theta_0$ (or $\lambda\geq 1)$, $R(\underline{\theta},\underline{\delta}_{\alpha}^{M})$ is minimized at
	\begin{align}\label{eq:5.5}
		\alpha=	\alpha(\lambda)&=\frac{p_1}{p_1+p_2}-\frac{p_2}{p_1+p_2} \frac{(\lambda-1) \,E_{\underline{\theta}}[Z_1(c_{0,1}-c_{0,2} \lambda Z) I_{(0,\frac{c_{0,1}}{c_{0,2}\lambda})}(Z)]}{E_{\underline{\theta}}[Z_1^2 (c_{0,1}-c_{0,2} \lambda Z)^2   \nonumber I_{(0,\frac{c_{0,1}}{c_{0,2}\lambda})}(Z)]}\\
		&=\frac{p_1}{p_1+p_2}-\frac{p_2}{p_1+p_2} \alpha_1(\lambda), \;\;\; (say),
	\end{align}
	where, for $h_1(z)=E[Z_1\vert Z=z]$ and $h_2(z)=E[Z_1^{^2}\vert Z=z]$,
	\begin{align*}
		\alpha_1(\lambda)&=\frac{(\lambda-1) \,E_{\underline{\theta}}[h_1(Z)(c_{0,1}-c_{0,2} \lambda Z) I_{(0,\frac{c_{0,1}}{c_{0,2}\lambda})}(Z)]}{E_{\underline{\theta}}[h_2(Z) (c_{0,1}-c_{0,2} \lambda Z)^2   \nonumber I_{(0,\frac{c_{0,1}}{c_{0,2}\lambda})}(Z)]} \qquad \qquad\\ \nonumber
		&= \frac{(\lambda-1) \,\int_{0}^{\frac{c_{0,1}}{c_{0,2}\lambda}} h_1(z)\, (c_{0,1}-c_{0,2}\lambda z)\,f_Z\left(z\right) dz}{\int_{0}^{\frac{c_{0,1}}{c_{0,2}\lambda}}h_2\left(z\right) (c_{0,1}-c_{0,2}\lambda z)^2 \,f_Z\left(z\right) dz}\\
		&= \frac{(\lambda-1) \,\int_{0}^{1} h_1\left(\frac{c_{0,1} t}{c_{0,2}\lambda}\right)\, (1-t)\,f_Z\left(\frac{c_{0,1} t}{c_{0,2}\lambda}\right) dt}{c_{0,1}\int_{0}^{1}h_2\left(\frac{c_{0,1} t}{c_{0,2}\lambda}\right) (1-t)^2 \,f_Z\left(\frac{c_{0,1} t}{c_{0,2}\lambda}\right) dt}\\
		&=\frac{1}{c_{0,1}}E_{\lambda}[k(S_{\lambda},\lambda)],\;\; \lambda\geq 1;
	\end{align*}
	here $k(t,\lambda)=\frac{(\lambda-1)h_1\left(\frac{c_{0,1} t}{c_{0,2}\lambda}\right)}{(1-t)\,h_2\left(\frac{c_{0,1} t}{c_{0,2}\lambda}\right)},\; 0<t<1,\; \lambda\geq 1,$ and $S_{\lambda}$ is a r.v. having the p.d.f.
	\begin{equation}\label{eq:5.6}
		h_{\lambda}(z)=\begin{cases}
			\frac{(1-t)^2 h_2\left(\frac{c_{0,1} t}{c_{0,2}\lambda}\right)f_Z\left(\frac{c_{0,1} t}{c_{0,2}\lambda}\right)}{\int_{0}^{1} (1-s)^2 h_2\left(\frac{c_{0,1} s}{c_{0,2}\lambda}\right) f_Z\left(\frac{c_{0,1} s}{c_{0,2}\lambda}\right)\,ds},& \text{if}\;\; 0<t<1\\
			0, &\text{otherwise}
		\end{cases},\;\;\lambda\geq 1.
	\end{equation}
	
	Clearly, if, for each $\theta\in(0,1)$, $\frac{h_2(\theta z)f_Z(\theta z)}{h_2(z)f_Z( z)}$ is increasing (decreasing) in $z\in \left(0,\frac{c_{0,1}}{c_{0,2}}\right)$, then $S_{\lambda_1}\leq_{lr}S_{\lambda_2}$ ($S_{\lambda_2}\leq_{lr}S_{\lambda_1}$) and, consequently, $S_{\lambda_1}\leq_{st}S_{\lambda_2}$ ($S_{\lambda_2}\leq_{st} S_{\lambda_1}$), whenever $1\leq \lambda_1<\lambda_2<\infty$. \vspace*{2mm}
	
	\noindent
	The following lemma, whose proof is immediate from Proposition 2.1, will be useful in proving the main result of this section.
	\\~\\ \textbf{Lemma 5.1.(a)} Suppose that, for each $\theta\in(0,1)$, $\frac{h_2(\theta z)f_Z(\theta z)}{h_2(z)f_Z( z)}$ is increasing (decreasing) in $z\in \left(0,\frac{c_{0,1}}{c_{0,2}}\right)$. Also, suppose that for every fixed $z\in (0,1)$, $k(z,\lambda)$ is increasing in $\lambda\in [1,\infty)$ and, for every fixed $\lambda\geq 1$, $k(z,\lambda)$ is increasing (decreasing) in $z\in (0,1)$. Then $\alpha_1(\lambda)$ is an increasing function (and hence $\alpha(\lambda)$ is a decreasing function) of $\lambda\in[1,\infty)$,
	$\inf_{\lambda\geq 1} \alpha(\lambda) = \frac{p_1}{p_1+p_2}-\frac{p_2}{p_1+p_2}\lim_{\lambda \to \infty}\alpha_1(\lambda)=\alpha_{\infty}, \text{ say, and }
	\sup_{\lambda\geq 1} \alpha(\lambda) =\frac{p_1}{p_1+p_2}-\frac{p_2}{p_1+p_2}\alpha_1(1)=\frac{p_1}{p_1+p_2}.$
	\\~\\ \textbf{(b)} Suppose that, for each $\theta\in(0,1)$, $\frac{h_2(\theta z)f_Z(\theta z)}{h_2(z)f_Z( z)}$ is increasing (decreasing) in $z\in \left(0,\frac{c_{0,1}}{c_{0,2}}\right)$. Also, suppose that, for every fixed $z\in (0,1)$, $k(z,\lambda)$ is decreasing in $\lambda\in [1,\infty)$ and, for every fixed $\lambda\geq 1$, $k(z,\lambda)$ is decreasing (increasing) in $z\in (0,1)$. Then $\alpha_1(\lambda)$ is a decreasing function (and hence $\alpha(\lambda)$ is an increasing function) of $\lambda\in[1,\infty)$,
	$\inf_{\lambda\geq 1} \alpha(\lambda) = \frac{p_1}{p_1+p_2}-\frac{p_2}{p_1+p_2}\alpha_1(1)=\frac{p_1}{p_1+p_2}, \text{ say, and }
	\sup_{\lambda\geq 1} \alpha(\lambda) =\frac{p_1}{p_1+p_2}-\frac{p_2}{p_1+p_2}\lim_{\lambda \to \infty}\alpha_1(\lambda)=\alpha_{\infty}.$
	\vspace*{2mm}
	
	Now we present the main result of this section. We use the convention that $[\frac{p_1}{p_1+p_2},\alpha_{\infty}]=[\frac{p_1}{p_1+p_2},\infty),$ if $\alpha_{\infty}=\infty$. The proof of the theorem, being similar to that of Theorem 4.1, is being omitted.
	\vspace*{2mm}

	\noindent  \textbf{Theorem 5.1. (a)}  Suppose that assumptions of Lemma 5.1 (a) hold. Then the estimators that are admissible within the class $\mathcal{D}=\{\underline{\delta}_{\alpha}^{M}=(\delta_{1,\alpha}^{M},\delta_{2,\alpha}^{M}): \alpha \in \Re\}$ are $\{\underline{\delta}_{\alpha}^{M}=(\delta_{1,\alpha}^{M},\delta_{2,\alpha}^{M}): \alpha \in [\alpha_{\infty},\frac{p_1}{p_1+p_2}]\}$. Moreover, for $-\infty<\alpha_1<\alpha_2\leq \alpha_{\infty}$ or $\frac{p_1}{p_1+p_2} \leq \alpha_2 < \alpha_1<\infty$, the estimator $\underline{\delta}_{\alpha_2}^{M}(\underline{X})$ dominates the estimator $\underline{\delta}_{\alpha_1}^{M}(\underline{X})$, for any $\underline{\theta}\in \Theta_0$.
	\vspace*{2mm}
	
	\noindent  \textbf{(b)}  Suppose that assumptions of Lemma 5.1 (b) hold. Then the estimators that are admissible within the class $\mathcal{D}=\{\underline{\delta}_{\alpha}^{M}=(\delta_{1,\alpha}^{M},\delta_{2,\alpha}^{M}): \alpha \in \Re\}$ are $\{\underline{\delta}_{\alpha}^{M}=(\delta_{1,\alpha}^{M},\delta_{2,\alpha}^{M}): \alpha \in [\frac{p_1}{p_1+p_2},\alpha_{\infty}]\}$. Moreover, for $-\infty<\alpha_1<\alpha_2\leq \frac{p_1}{p_1+p_2} $ or $\alpha_{\infty} \leq \alpha_2 < \alpha_1<\infty$, the estimator $\underline{\delta}_{\alpha_2}^{M}(\underline{X})$ dominates the estimator $\underline{\delta}_{\alpha_1}^{M}(\underline{X})$, for any $\underline{\theta}\in \Theta_0$.

	\subsection{Applications}
	\label{sec:3.3}
	\noindent
	
	\vspace*{3mm}

	\noindent
	Now we illustrate some applications of Theorem 5.1.
	
	\vspace*{3mm}
	
	\noindent \textbf{Example 5.1.1.} Let $\underline{X}=(X_1,X_2)$ be a random vector having the p.d.f. \eqref{eq:5.1}, with $f(z_1,z_2)=f_1(z_1) f_2(z_2),\;\; \underline{z}=(z_1,z_2)\in \Re^2$ and
	$$f_i(z)=\begin{cases}  \frac{e^{-z} z^{\alpha_i-1}}{\Gamma{\alpha_i}}, &\text{ if} \;\; z>0 \\
		0, & \text{ otherwise}  \end{cases},\;\; i=1,2.  $$
	Here $\alpha_1>0$ and $\alpha_2>0$ are known shape parameters. We have $c_{0,i}=\frac{E[Z_i]}{E[Z_i^2]}=\frac{1}{\alpha_i+1},\;i=1,2,$ and the BSEE of $\theta_i$ is $\delta_{i,0}(\underline{X})=\frac{X_i}{\alpha_i+1},\;i=1,2$. The p.d.f. of $Z=\frac{Z_2}{Z_1}$ is 
	
	$$f_Z(z)= \begin{cases} \frac{\Gamma{(\alpha_1+\alpha_2)}}{\Gamma{\alpha_1} \Gamma{\alpha_2}}\frac{z^{\alpha_2-1}}{(1+z)^{\alpha_1+\alpha_2}}, &\text{ if} \;\; 0<z<\infty \\
		0, & \text{ otherwise}  \end{cases},  $$
	and the conditional p.d.f. of $Z_1$ given $Z=z$ ($0<z<\infty$) is
	$$ f_{Z_1|Z}(z_1|z)=\begin{cases}
		\frac{ (1+z)^{\alpha_1+\alpha_2}\, z_1^{\alpha_1+\alpha_2-1}\,  e^{-(1+z)z_1}}{\Gamma{(\alpha_1+\alpha_2)}}, &\text{ if} \;\; 0<z_1<\infty \\
		0, & \text{ otherwise}  \end{cases}. $$
	Clearly, $h_1(z)=E[Z_1|Z=z]=\frac{\alpha_1+\alpha_2}{1+z},\; z>0,$
	$h_2(z)=E[Z_1^2|Z=z]=\frac{(\alpha_1+\alpha_2+1)(\alpha_1+\alpha_2)}{(1+z)^2},\; z>0,$ and, for $0<\theta<1$, $\frac{h_2(\theta z)f_Z(\theta z)}{h_2(z) f_Z(z)}$ is increasing in $z\in \left(0,\frac{c_{0,1}}{c_{0,2}}\right]$. For any fixed $\lambda\geq 1$,
	$$ k(z,\lambda)=\frac{(\lambda-1)h_1\left(\frac{c_{0,1} z}{c_{0,2}\lambda}\right)}{(1-z)\,h_2\left(\frac{c_{0,1} z}{c_{0,2}\lambda}\right)}=\frac{(\lambda-1)}{(1-z)(\alpha_1+\alpha_2+1)}\left(1+\frac{c_{0,1} z}{c_{0,2}\lambda}\right)$$
	is increasing in $z\in (0,1)$ and, for any fixed $z\in (0,1)$, $k(z,\lambda)$ is increasing in $\lambda\geq 1$. Using Lemma 5.1 (a), we have
	\begin{align*}
		\sup_{\lambda\geq 1}\alpha(\lambda)&=\frac{p_1}{p_1+p_2},\\~\\
		\text{and}\quad 
		\inf_{\lambda\geq 1}\alpha(\lambda)&=\frac{p_1}{p_1+p_2}-\frac{p_2}{p_1+p_2}\lim_{\lambda \to \infty}\alpha_1(\lambda)\\
		&=\frac{p_1}{p_1+p_2}-\frac{p_2}{p_1+p_2}\lim_{\lambda \to \infty}\frac{(\lambda-1) \,\int_{0}^{1} h_1\left(\frac{c_{0,1} t}{c_{0,2}\lambda}\right)\, (1-t)\,f_Z\left(\frac{c_{0,1} t}{c_{0,2}\lambda}\right) dt}{c_{0,1}\int_{0}^{1}h_2\left(\frac{c_{0,1} t}{c_{0,2}\lambda}\right) (1-t)^2 \,f_Z\left(\frac{c_{0,1} t}{c_{0,2}\lambda}\right) dt}\\
		&=-\infty.
	\end{align*}
	Using Theorem 5.1 (a), it follows that the estimators $\{\underline{\delta}_{\alpha}^{M}:\alpha\in (-\infty,\frac{p_1}{p_1+p_2}]\}$ are admissible within the class $\{\underline{\delta}_{\alpha}^{M}:-\infty<\alpha<\infty\}$ of isotonic regression estimators of $(\theta_1,\theta_2)$. Moreover the estimators $\{\underline{\delta}_{\alpha}^{M}:\alpha>\frac{p_1}{p_1+p_2}\}$ are inadmissible and, for $\frac{p_1}{p_1+p_2}\leq \alpha_2<\alpha_1$, the estimator $\underline{\delta}_{\alpha_2}^{M}$ dominates the estimator $\underline{\delta}_{\alpha_1}^{M}$. Using Theorem 3.1, we conclude that estimators $\{\underline{\delta}_{\alpha}^{M}:\max\{0,\frac{p_1-p_2}{p_1}\}\leq \alpha<1\}$ dominate the BSEE $(\frac{1}{\alpha_1+1}X_1,\frac{1}{\alpha_2+1}X_2)$.
	
	\vspace*{3mm}

	\noindent \textbf{Example 5.1.2.} Let $\underline{X}=(X_1,X_2)$ be a random vector with p.d.f. \eqref{eq:5.1}, where $f(z_1,z_2)=f_1(z_1) f_2(z_2),\; \underline{z}=(z_1,z_2)\in \Re_{++}^2$ and, for positive constants $\alpha_1$ and $\alpha_2$,
	$$f_i(z)=\begin{cases}  \alpha_i z^{\alpha_i-1}, &\text{ if} \;\; 0<z<1 \\
		0, & \text{ otherwise}  \end{cases},\;\; i=1,2.  $$
	Here $c_{0,i}=\frac{E[Z_i]}{E[Z_i^2]}=\frac{\alpha_i+2}{\alpha_i+1},\;i=1,2,$ and the BSEE of $\theta_i$ is $\delta_{i,0}(\underline{X})=\frac{\alpha_1+2}{\alpha_i+1}X_i,\;i=1,2$. The p.d.f. of $Z=\frac{Z_2}{Z_1}$ is 
	$$f_Z(z)= \frac{\alpha_1 \alpha_2}{\alpha_1+\alpha_2} z^{\alpha_2-1} \left(\min\bigg\{1,\frac{1}{z}\bigg\}\right)^{\alpha_1+\alpha_2},\;\; z>0,  $$
	the conditional p.d.f. of $Z_1$ given $Z=z$ ($0<z<\infty$) is
	$$ f_{Z_1|Z}(z_1|z)=\begin{cases}
		\frac{(\alpha_1+\alpha_2) \,z_1^{\alpha_1+\alpha_2-1}}{ \left(\min\big\{1,\frac{1}{z}\big\}\right)^{\alpha_1+\alpha_2}}, &\text{ if} \;\; 0<z_1<\min\big\{1,\frac{1}{z}\big\} \\
		0, & \text{ otherwise}  \end{cases} , $$ 
	$h_1(z)=E[Z_1|Z=z]=\frac{\alpha_1+\alpha_2}{\alpha_1+\alpha_2+1}\min\big\{1,\frac{1}{z}\big\},\, z>0, \text{ and } 
	h_2(z)=E[Z_1^2|Z=z]=\frac{\alpha_1+\alpha_2}{\alpha_1+\alpha_2+2}(\min\big\{1,\frac{1}{z}\big\})^2,$ $\, z>0.$ Clearly, for $0<\theta<1$, $\frac{h_2(\theta z)f_Z(\theta z)}{h_2(z) f_Z(z)}$ is increasing in $z\in \left(0,\frac{c_{0,1}}{c_{0,2}}\right]$. For any fixed $\lambda\geq 1$, 
	$$ k(z,\lambda)=\frac{(\lambda-1)h_1\left(\frac{c_{0,1} z}{c_{0,2}\lambda}\right)}{(1-z)\,h_2\left(\frac{c_{0,1} z}{c_{0,2}\lambda}\right)}=\frac{(\lambda-1)(\alpha_1+\alpha_1+2)}{(1-z)(\alpha_1+\alpha_2+1)}\max\bigg\{1,\frac{c_{0,1}z}{c_{0,2}\lambda}\bigg\}$$
	is increasing in $z$ and, for any fixed $z\in (0,1)$, $k(z,\lambda)$ is increasing in $\lambda\geq 1$. Using Lemma 5.1 (a), we get
	\begin{align*}
		\sup_{\lambda\geq 1}\alpha(\lambda)=\frac{p_1}{p_1+p_2},\qquad
		\text{and}\qquad 
		\inf_{\lambda\geq 1}\alpha(\lambda)=-\infty.
	\end{align*}
	
	Now, using Theorem 5.1 (a), it follows that the estimators $\{\underline{\delta}_{\alpha}^{M}:\alpha\in (-\infty,\frac{p_1}{p_1+p_2}]\}$ are admissible within the class $\{\underline{\delta}_{\alpha}^{M}:-\infty<\alpha<\infty\}$ of isotonic regression estimators of $(\theta_1,\theta_2)$. Moreover the estimators $\{\underline{\delta}_{\alpha}^{M}:\alpha>\frac{p_1}{p_1+p_2}\}$ are inadmissible and, for $\frac{p_1}{p_1+p_2}\leq \alpha_2<\alpha_1$, the estimator $\underline{\delta}_{\alpha_2}^{M}$ dominates the estimator $\underline{\delta}_{\alpha_1}^{M}$. Also, using Theorem 3.1, we conclude that the class of estimators $\{\underline{\delta}_{\alpha}^{M}:\max\{0,\frac{p_1-p_2}{p_1}\}\leq \alpha<1\}$ dominate the BSEE $(\frac{\alpha_1+2}{\alpha_1+1}X_1,\frac{\alpha_2+2}{\alpha_2+1}X_2)$.

	\subsection{\textbf{Simulation Study For Estimation of Scale Parameters $(\theta_1,\theta_2)$}} 
	\label{sec:4.4}
	\noindent

	\vspace*{3mm}

	In Example 5.1.1, we have considered two independent gamma distributions with unknown order restricted scale parameters (i.e., $\theta_1\leq \theta_2$) and known shape parameters ($\alpha_1>0$ and $\alpha_2>0$). For simultaneous estimation of scale parameters ($\theta_1$,$\theta_2$), under the sum of the squared error loss functions \eqref{eq:5.2},
	we have shown that the isotonic regression estimators (IREs) $\{\underline{\delta}_{\alpha}^{M}: \alpha\in(-\infty,\frac{p_1}{p_1+p_2}] \}$ are admissible within the class $\mathcal{D}=\{\underline{\delta}_{\alpha}^{M}: -\infty<\alpha<\infty \}$. Also, the estimators $\{\underline{\delta}_{\alpha}^{M}: \max\{0,\frac{p_1-p_2}{p_1}\}\leq \alpha<1 \}$ dominate the BSEE $\underline{\delta}_{0}(\underline{X})=(c_{0,1}X_1,c_{0,2}X_2)$, where $c_{0,i}=\frac{1}{\alpha_i+1},\;i=1,2$. In this section, for simplicity, we take $p_1=1$ and $p_2=1$. To further evaluate the performances of various estimators under the loss function 	$L(\underline{\theta},\underline{a})=(a_1-\theta_1)^2 + (a_2-\theta_2)^2, \;\underline{\theta}=(\theta_1,\theta_2)\in\Theta_0,\;\underline{a}=(a_1,a_2)\in\Re_{++}^2$, in this section, we compare the risk performances of isotonic regression estimators (IREs) $\underline{\delta}_{1}^{M}(\underline{X})=(\delta_{1,1}^{M}(\underline{X}),\delta_{2,1}^{M}(\underline{X}))=(c_{0,1}X_1,c_{0,2}X_2)$, $\underline{\delta}_{0.75}^{M}(\underline{X})=(\delta_{1,0.75}^{M}(\underline{X}),\delta_{2,0.75}^{M}(\underline{X}))$, $\underline{\delta}_{0.5}^{M}(\underline{X})=(\delta_{1,0.5}^{M}(\underline{X}),\delta_{2,0.5}^{M}(\underline{X}))$, $\underline{\delta}_{0}^{M}(\underline{X})=(\delta_{1,0}^{M}(\underline{X}),\delta_{2,0}^{M}(\underline{X}))$ and the restricted MLE $\underline{\delta}_R(\underline{X})=\left(\min\Big\{\frac{X_1}{\alpha_1},\frac{X_1+X_2}{\alpha_1+\alpha_2}\Big\},\max\Big\{\frac{X_2}{\alpha_2},\frac{X_1+X_2}{\alpha_1+\alpha_2}\Big\}\right)$, numerically, through the Monte Carlo simulations. 
	\vspace*{2mm}
	
	For simulations, we generated 50000 samples of size 1 each from relevant gamma distributions and computed the simulated risks of estimators $\underline{\delta}_{1}^{M}(\underline{X})=(X_1-\sigma_1,X_2-\sigma_2)$, $\underline{\delta}_{0.75}^{M}(\underline{X})$, $\underline{\delta}_{0.5}^{M}(\underline{X})$, $\underline{\delta}_{0}^{M}(\underline{X})$ and $\underline{\delta}_{R}^{M}(\underline{X})$ for different values of shape parameters $(\alpha_1,\alpha_2)$.

	\FloatBarrier
	\begin{figure}[h!]
		\begin{subfigure}{0.48\textwidth}
			\centering
			\includegraphics[width=85mm,scale=1.2]{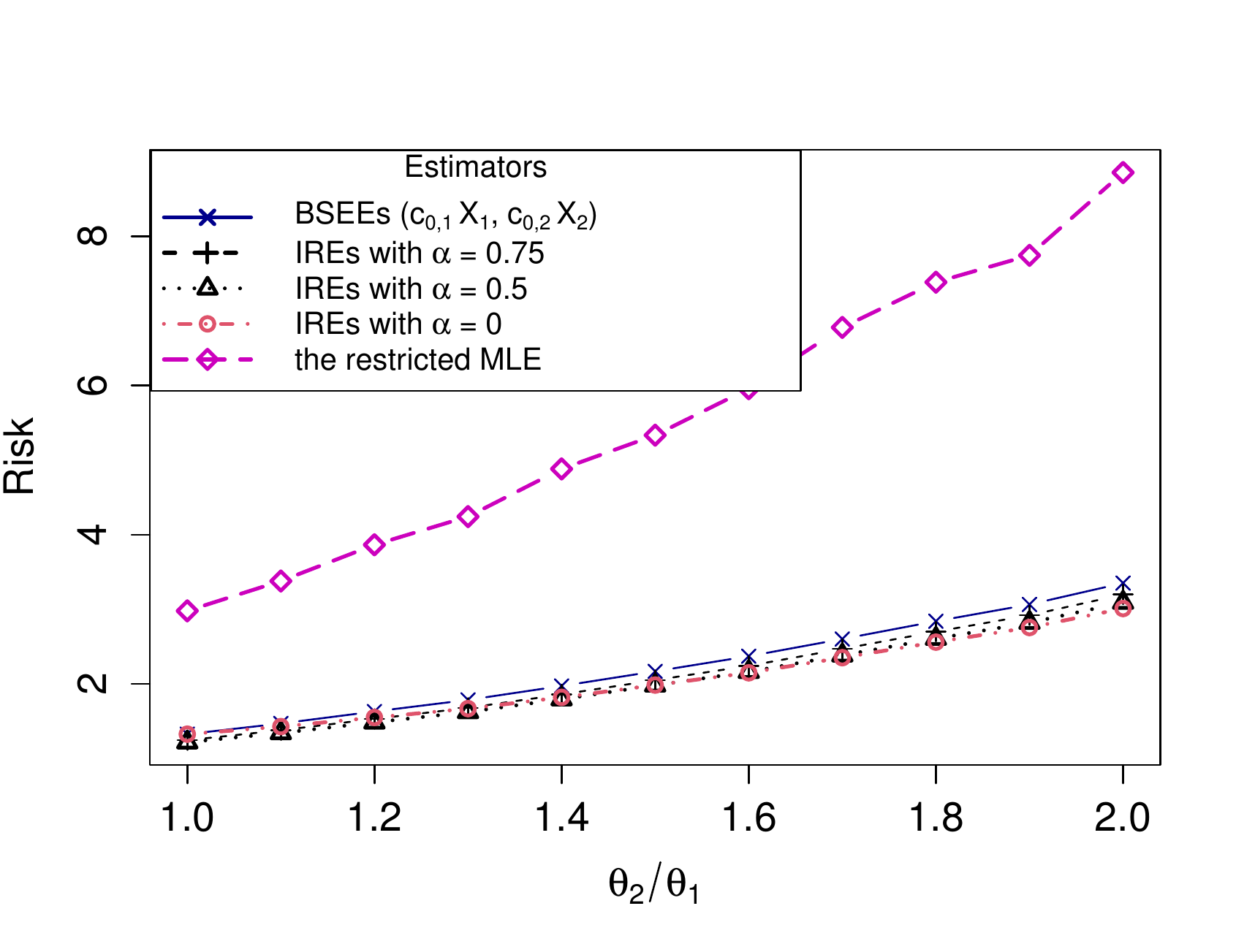} 
			\caption{$\alpha_1=0.5$ and $\alpha_2=0.5$.} 
			\label{fig1:a} 
		\end{subfigure}
		\begin{subfigure}{0.48\textwidth}
			\centering
			\includegraphics[width=85mm,scale=1.2]{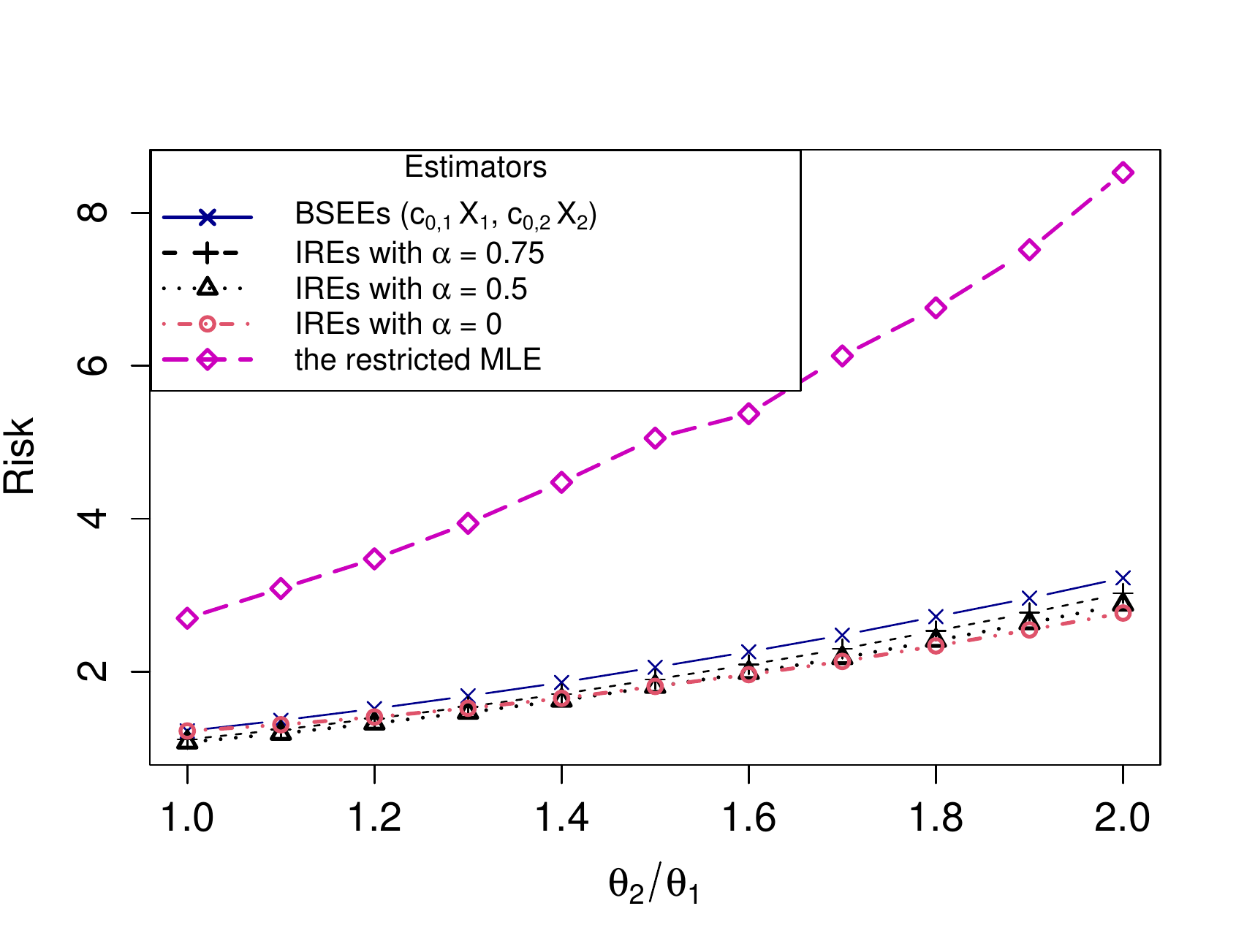} 
			\caption{$\alpha_1=0.8$ and $\alpha_2=0.5$.} 
			\label{fig1:b} 
		\end{subfigure}
		\\	\begin{subfigure}{0.48\textwidth}
			\centering
			
			\includegraphics[width=85mm,scale=1.2]{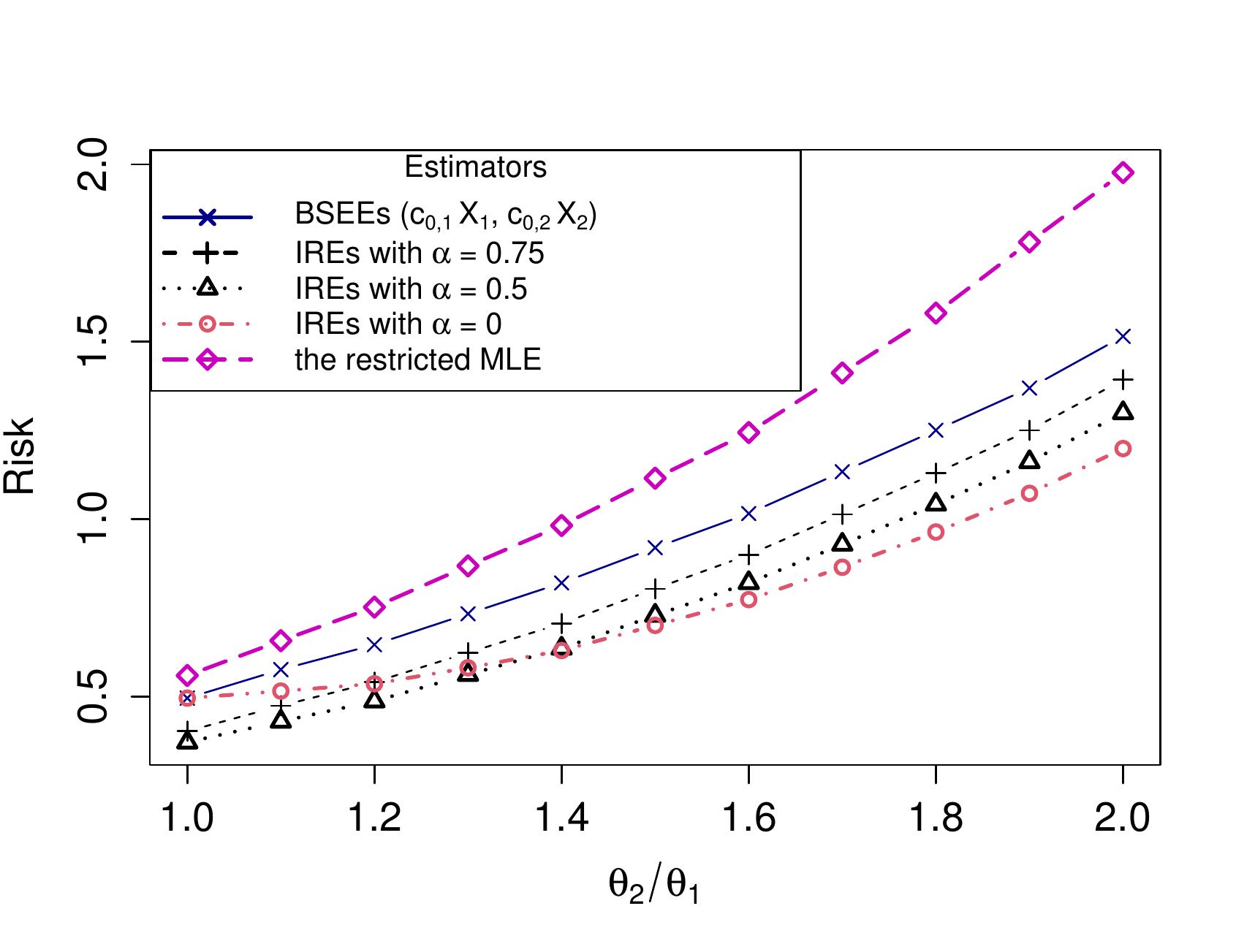} 
			
			\caption{$\alpha_1=5$ and $\alpha_2=2$.} 
			\label{fig1:c} 
		\end{subfigure}
		\begin{subfigure}{0.48\textwidth}
			\centering
			\includegraphics[width=85mm,scale=1.2]{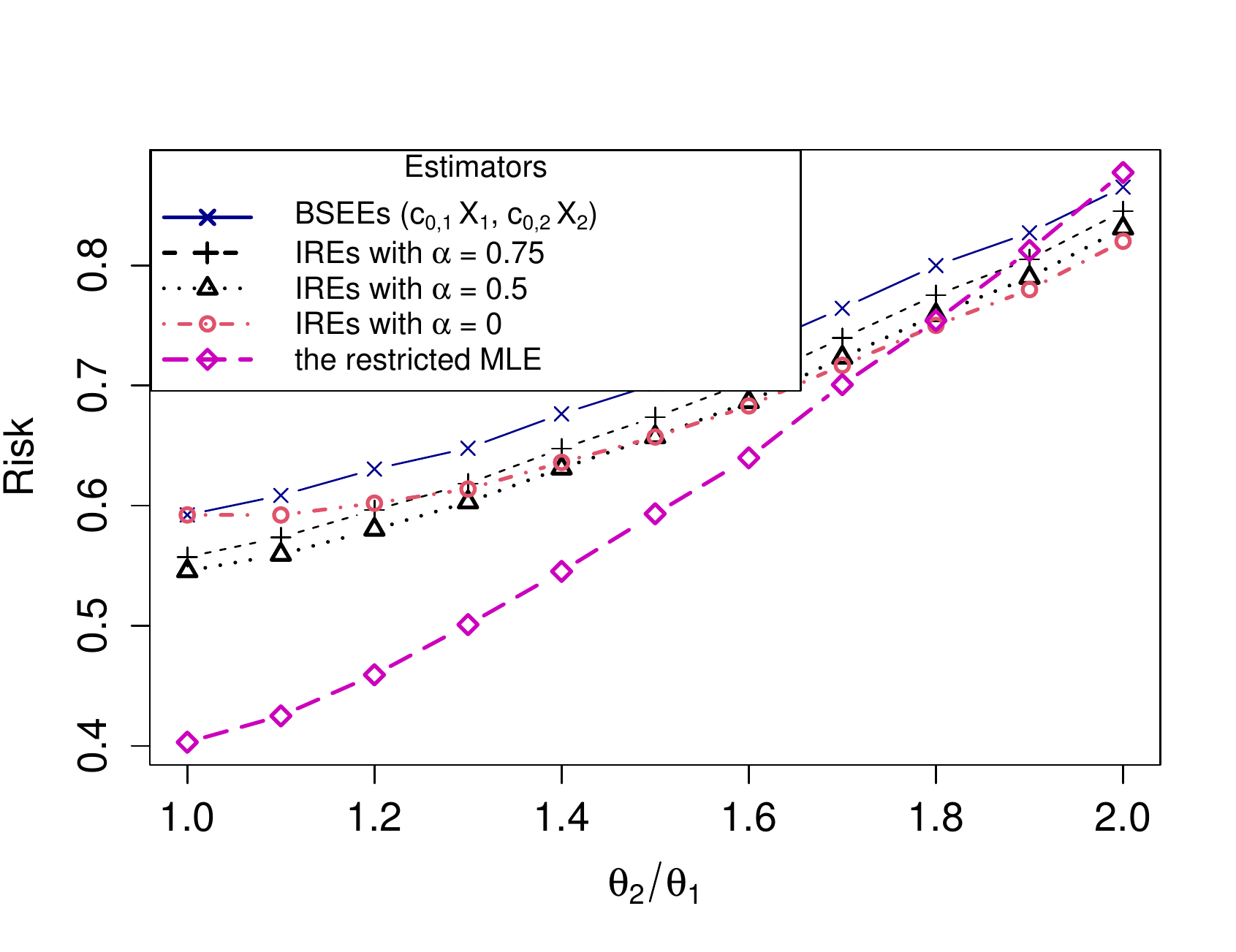} 
			\caption{ $\alpha_1=1$ and $\alpha_2=10$.} 
			\label{fig1:d}  
		\end{subfigure}
		\\	\begin{subfigure}{0.48\textwidth}
			\centering
			
			\includegraphics[width=85mm,scale=1.2]{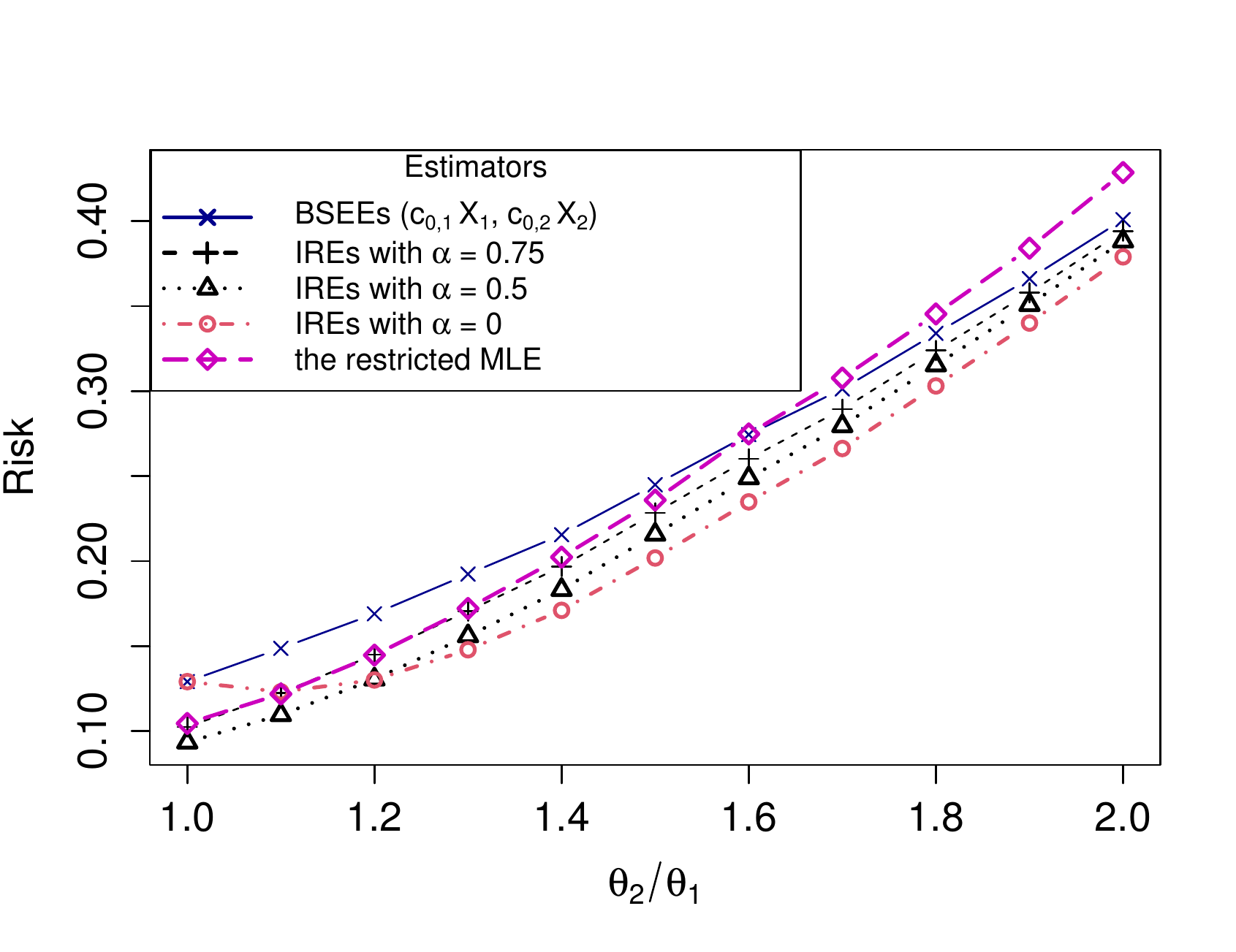} 
			
			\caption{$\alpha_1=25$ and $\alpha_2=10$.} 
			\label{fig1:e} 
		\end{subfigure}
		\begin{subfigure}{0.48\textwidth}
			\centering
			\includegraphics[width=85mm,scale=1.2]{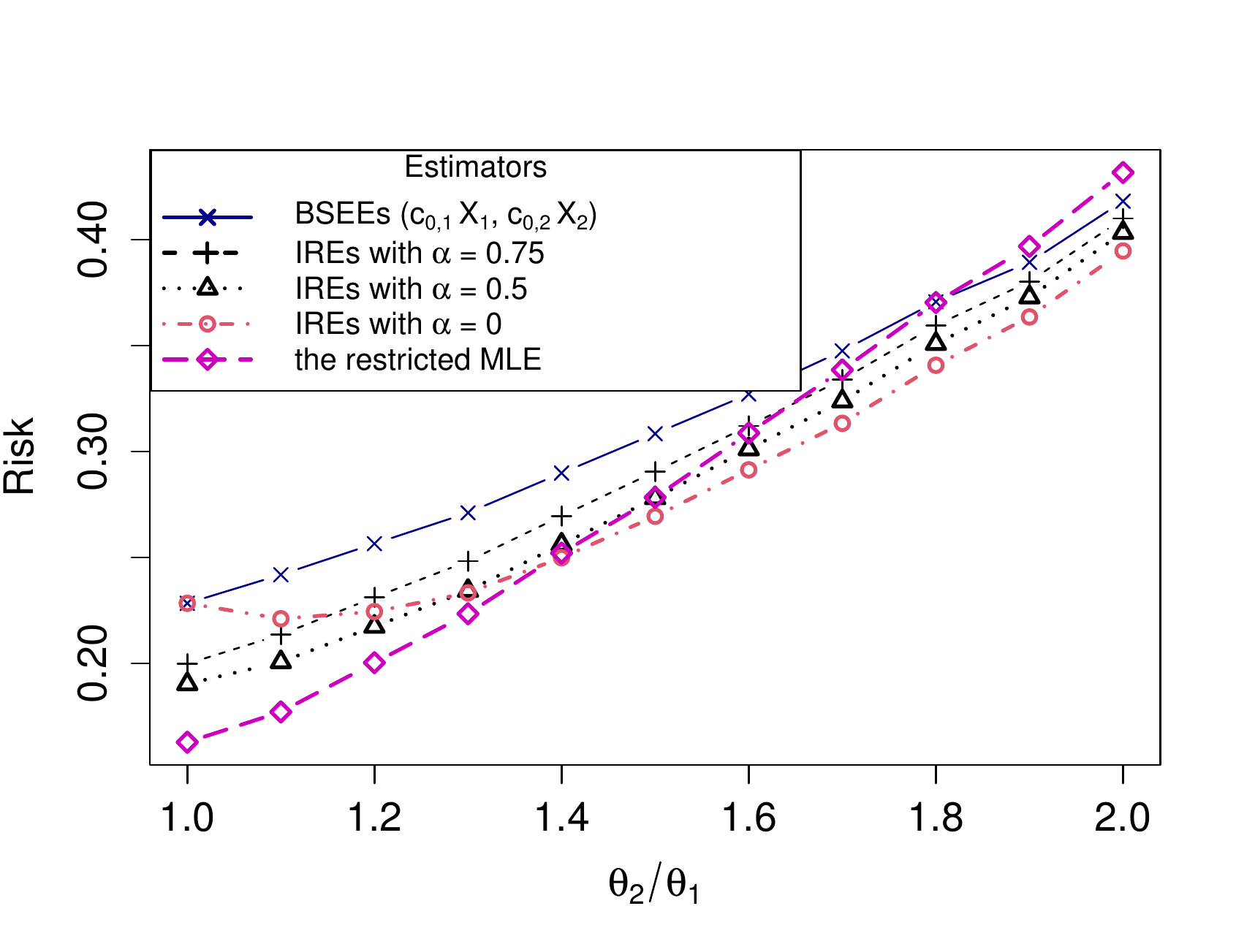} 
			\caption{$\alpha_1=5$ and $\alpha_2=10$.} 
			\label{fig1:f}  
		\end{subfigure}
		\caption{Risk plots of estimators $\underline{\delta}_{1}^{M}(\underline{X})=(c_{0,1}X_1,c_{0,2}X_2)$, $\underline{\delta}_{0.75}^{M}(\underline{X})$, $\underline{\delta}_{0.5}^{M}(\underline{X})$, $\underline{\delta}_{0}^{M}(\underline{X})$ and $\underline{\delta}_{R}^{M}(\underline{X})$ against the values of $\lambda=\theta_2/\theta_1$.}
		\label{fig2}
	\end{figure}
	\FloatBarrier

	The simulated values of risks of various estimators are plotted in Figure 3. The following observations are evident from Figure 3:\vspace*{2mm}

	\noindent (i) The risk function values of estimators $\underline{\delta}_{0.75}^{M}(\underline{X})$, $\underline{\delta}_{0.5}^{M}(\underline{X})$ and $\underline{\delta}_{0}^{M}(\underline{X})$ are less than the risk function values of the BSEE $(X_1-\sigma_1,X_2-\sigma_2)$, which is in 
	conformity with theoretical findings of Example 5.1.1 and Theorem 5.1.
	\noindent \\~\\(ii) From Figure 3, we can observe that the IREs $\underline{\delta}_{0.5}^{M}(\underline{X})$, $\underline{\delta}_{0}^{M}(\underline{X})$ and $\underline{\delta}_{R}^{M}(\underline{X})$ are not comparable and the IREs $\underline{\delta}_{1}^{M}(\underline{X})$ and $\underline{\delta}_{0.75}^{M}(\underline{X})$ are inadmissible. This is in conformity with theoretical findings of Example 5.1.1.
	\\~\\(iii) When values of $\alpha_1$ and $\alpha_2$ are large and $\alpha_1<<\alpha_2$, for small values of $\theta_2/\theta_1$, the restricted MLE outperforms the other estimators. 
	\\~\\(iv) There is no clear cut winner between various estimators, but performance of the BLEE $\underline{\delta}_{1}^{M}(\underline{X})=(X_1-\sigma_1,X_2-\sigma_2)$ and $\underline{\delta}_{0.75}^{M}(\underline{X})$ are worse than other estimators. Also, when $\alpha_1\leq 1$ and $\alpha_2\leq 1$, the performance of the restricted MLE is the worst among other estimators.

	\section*{\textbf{Disclosure statement}}
	
	There is no conflict of interest by authors.

	\section*{\textbf{Funding}}
	
	This work was supported by the [Council of Scientific and Industrial Research (CSIR)] under Grant [number 09/092(0986)/2018].

	\bibliographystyle{apalike}
	\bibliography{Paper5}		
	
\end{document}